\documentclass[11pt]{article}
\usepackage{amsmath, amsthm, amssymb, eucal, stmaryrd, setspace}
\usepackage[all]{xy}

\usepackage{url}
\numberwithin{equation}{subsection}
\swapnumbers
\newtheorem{theor}[equation]{Theorem}
\newtheorem{cor}[equation]{Corollary}
\newtheorem{lem}[equation]{Lemma}
\newtheorem{pro}[equation]{Proposition}
\newtheorem{point}[equation]{}
\theoremstyle{remark}
\newtheorem{defn}[equation]{Definition}

\newtheorem{rem}[equation]{Remark}
\newtheorem{gener}[equation]{Generalization}
\newtheorem{exa}[equation]{Example}

\newtheorem{notn}[equation]{Notation}

\newcommand{\into}{\hookrightarrow}
\newcommand{\Af}{\mathbb A}
\newcommand{\CC}{\mathbb C}

\newcommand{\GG}{\mathbb G}

\newcommand{\HH}{\mathbb H}

\newcommand{\NN}{\mathbb N}
\newcommand{\PP}{\mathbb P}
\newcommand{\QQ}{\mathbb Q}

\newcommand{\ZZ}{\mathbb Z}
\newcommand{\Ocal}{\mathcal O}
\newcommand{\Ccal}{\mathcal C}

\newcommand{\Ecal}{\mathcal E}

\newcommand{\isoto}{\xrightarrow{\ \sim\ } }
\newcommand{\rest}[1]{{}{\textstyle|}_{#1}}
\newcommand{\size}[1]{{\#(}{#1}{)}}
\newcommand{\id}{\operatorname{id}}
\newcommand{\Aut}{\operatorname{Aut}}

\newcommand{\hcf}{\operatorname{hcf}}

\newcommand{\Isom}{\sta {Isom}}

\newcommand{\mult}{\operatorname{mult}}

\newcommand{\Pic}{\operatorname{Pic}}

\newcommand{\Spec}{\operatorname{Spec}}

\newcommand{\al}{\alpha}
\newcommand{\be}{\beta}

\newcommand{\ep}{\varepsilon}
\newcommand{\fie}{\varphi}

\newcommand{\la}{\lambda}

\newcommand{\coker}{\operatorname{coker}}
\newcommand{\sh}{\operatorname{sh}}

\providecommand{\abs}[1]{\lvert#1\rvert}
\newcommand{\coa}{\abs}
\newcommand{\sta}{\mathsf}
\newcommand{\stC}{\sta C}
\newcommand{\stL}{\sta L}

\def\pmmu{{\pmb \mu}}
\def\sm{^{\rm sm}}

\def\ol{\overline}
\newcommand{\lvec}[1]{\vec{#1}\,}

\def\wt{\widetilde}

\newcommand{\BBB}{\sta B}
\newcommand{\MMM}{{\sta M}}

\newcommand{\MMMbar}{\ol{\MMM}{}}

\begin{document}
\title{
\textbf
{Stable twisted curves and their $r$-spin structures}}
\author{Alessandro Chiodo\thanks
{Financially supported by
the Marie Curie Intra-European
Fellowship within the 6th European Community Framework Programme,
MEIF-CT-2003-501940.}
}
\maketitle

\begin{quote}
\abstract{\noindent
The object of this paper is the notion of
$r$-spin structure: a line bundle
whose $r$th power is
isomorphic to the canonical bundle.
Over the moduli functor ${\mathsf{M}}_g$ of
smooth genus-$g$ curves, $r$-spin structures
form a finite torsor
under the group
of $r$-torsion line bundles.
Over the moduli functor $\overline{\mathsf{M}}_g$ of stable curves,
$r$-spin structures
form an \'etale stack, but the finiteness
and the torsor structure are lost.

In the present work, we show how
this bad picture can be definitely
improved simply by placing the problem
in the category of Abramovich and Vistoli's twisted curves.
First, we find that within such category
there exist several different compactifications of
${\mathsf{M}}_g$; each one corresponds to a different multiindex
$\vec{l}=(l_0,l_1,\dots)$ identifying a notion of stability: $\vec{l}$-stability.
Then, we determine the suitable choices of $\vec{l}$ for
which $r$-spin structures form a finite torsor
over the moduli
of $\vec{l}$-stable curves.}
\end{quote}

\section{Introduction}
For any integer $r\ge 2$, spin structures
of order $r$ are natural generalizations
of theta characteristics: on a space, they are given by
a line bundle $L$ and an
isomorphism $f\colon L^{\otimes r}\isoto \omega$.
In this paper we focus on their moduli functor.

For a fixed integer $r\ge 2$, we work  over
$\Spec\ZZ[1/r]$.

\subsection{Smooth curves: the $r$th roots form a torsor}
For $g\ge 2$ and $2g-2\in r\ZZ$,
the category of $r$-spin structures
on smooth genus-$g$ curves
forms a Deligne--Mumford stack
$\MMM_g^{\omega, r}$,
finite and \'etale on $\MMM_g$, which we
write as
$$\MMM_g^{\omega, r}=\{(C,L,f)\mid f\colon L^{\otimes r} \isoto \omega_C\}
_{/\cong}\longrightarrow \MMM_g.$$
In fact $\MMM_g^{\omega, r}$ is a finite torsor under
the finite group stack of
$r$-torsion line bundles on smooth
genus-$g$ curves:
$$\MMM_g^{\Ocal, r}=\{(C,L,f)\mid f\colon L^{\otimes r} \isoto \Ocal_C\}
_{/\cong}\longrightarrow \MMM_g.$$

\subsection{Stable curves: the torsor is lost}
When we extend the study of $r$-spin structures to
the category of stable curves $\MMMbar_g$,
the properness and the
torsor structure are lost.

First, consider the category
$\MMMbar_g^{\Ocal, r}$ of $r$-torsion line bundles
on stable curves.
As above, it forms an \'etale stack
on $\MMMbar_g$ which is equipped
with a group structure.
However, $\MMMbar_g^{\Ocal, r}\to \MMMbar_g$ is not proper.
Indeed, since $\MMMbar_g^{\Ocal, r}$ is \'etale
and the generic fibre contains $r^{2g}$ points,
one can check that the valuative criterion
fails by exhibiting an example of a geometric fibre
with less than $r^{2g}$ points.
\begin{exa}\label{exa:intro1}
Consider an irreducible curve of genus $g$ with only one node.
Note that the set of roots of $\Ocal$ consists of $r^{2g-1}$ elements.
Indeed, and more generally, for any stable curve $C$
the group of $r$-torsion line bundles $(\Pic C)_r$ fits in the exact sequence
\begin{equation}\label{eq:stable_sequence}
1\to \pmmu_r\to (\pmmu_r)^{\# V}\to (\pmmu_r)^{\# E}\to
(\Pic C)_r\to (\Pic C^{\nu})_r\to 1\end{equation}
where  $V$ and $E$ are the
sets of irreducible components and of singularities of
the curve $C$, whereas
$C^{\nu}$
is the normalization of $C$. So, we get
\begin{equation}\label{eq:intrortors}
\#{(\Pic C)_r}=r^{2g-1+\size{V}-\size{E}}.
\end{equation}
\end{exa}

Second, write
$\MMMbar_g^{\omega, r}$
for the category of $r$-spin structures on
stable curves. In fact, the morphism to $\MMMbar_g$ is \'etale;
however,
it is not proper and it is not
a torsor on $\MMMbar_g$.
Indeed, the following example shows that the
morphism to $\MMMbar_g$ is not surjective.
\begin{exa}
Let $C$ be a curve with only one node and
two irreducible components of genus $g-1$ and $1$, respectively.
Then there are no $r$th roots of $\omega$ on $C$.
This happens because the degree of $\omega$ is
$1$ on the genus-$1$ component of $C$.
Indeed, recall that the degree
of the dualizing sheaf $\omega$ on an irreducible
component $C'$ of genus $i$ is $2i-2+\size {N}$
where $N$ is the set of points where
$C'$ meets the rest of the curve
(by ``genus'' we always mean the \emph{arithmetic} genus, \S2.3).
\end{exa}

In the recent years, the interest in moduli of $r$-spin structures
has been
revived by Witten's conjecture \cite{Wi}, which relates certain
enumerative properties of $r$-spin structures   to the
Gelfand--Diki\u\i\  hierarchy. The conjecture is a generalization
of the Kontsevich--Witten Theorem \cite{Wi_2} \cite{Ko} and has been lately proven in \cite{FSZ}.
This result opens the way to further investigations of Gromov--Witten
$r$-spin theory.

The original formulation of the relevant enumerative properties
was only sketched by Witten in \cite{Wi}.
A rigourous definition requires---first of all---a
compactification of $\MMM_g^{\omega, r}$ (once a suitable
compactification is given, the numerical invariants can
be defined using \cite{PV} or \cite{Ch}, see Proposition \ref{pro:wtcc}).

In the existing literature,
there are several solutions to the
problem of compactifying $\MMM_g^{\omega, r}$: they consist
in enlarging the category $\MMM_g^{\omega, r}$
of smooth $r$-spin curves $(C,L,f)$
to a new category fibred over $\MMMbar_g$.
In \cite{Ja_1} and \cite{Ja_geom}, Jarvis allows noninvertible sheaves.
In \cite{Co} for $r=2$ and in \cite{CCC} for all $r$, Cornalba, Caporaso, and Casagrande
take as new objects line bundles on
semistable curves. In \cite{AJ}
Abramovich and Jarvis realize the same category as in \cite{Ja_geom} in
terms of stack-theoretic curves.
In all these compactifications
the torsor structure is lost over $\MMMbar_g$,
because ramification occurs at the new points
\cite[Thm.~2.4.2]{Ja_geom}
and \cite[\S3, \S4.1]{CCC}.

\subsection{Placing the problem in the context of twisted curves: $\lvec{l}$-stability}
We consider the category $\wt \MMM_g$
of twisted curves, which are, over an algebraically closed field,
stack-theoretic curves whose
smooth geometric locus is represented by a scheme
and whose stabilizers at the nodes have finite order
(see Abramovich and Vistoli \cite{AV} or
\S\ref{notn:twisted} for the definition over a
base scheme $X$).

Olsson  shows that
$\wt \MMM_g$ forms an algebraic stack, \cite{Ol}.
However, the stack $\wt \MMM_g$ is nonseparated. Indeed,
a twisted curve $\stC$ over a discrete valuation ring $R$
with smooth generic fibre $\stC_K$ is isomorphic to
its coarse space over the field of fractions
$K$ and may differ from it on the special fibre;
in this case  the coarse space ${\coa{\sta C}}$ and
the twisted curve $\stC$  are two nonisomorphic
twisted curves extending $\stC_K$ over $R$.
Therefore, the valuative criterion of separateness
fails.

We describe the condition of stability in the category
of twisted curves $\wt\MMM_g$.
For any multiindex
$\lvec{l}=(l_0,l_1,\dots,l_{\lfloor g/2\rfloor})$ of
invertible integers,
we say that a twisted curve is \emph{$\lvec{l}$-stable}
if its stabilizers
have order $l_i$ on nodes of type $i$ (the notion of type
of a node can be found in \cite{DM} and is recalled
in \eqref{notn:type}).
In this way, for each multiindex $\lvec{l}$,
we have a notion of stability, which corresponds to a new compactification
of $\MMM_g$ (the classical
Deligne--Mumford--Knudsen compactification
$\MMMbar_g$
corresponds to $\lvec{l}=(1,\dots,1)$).
In the following theorem,
we show that the compactifications obtained in this way
are all the compactifications of $\MMM_g$
inside $\wt\MMM_g$.

\vspace{0.2cm}
\noindent{\bf 
\ref{thm:olsson} Theorem.}
{\it
Let us denote by $\MMM_g(\lvec l)$ the category
of $\lvec{l}$-stable curves. It is contained in
$\wt \MMM_g$ and it contains $\MMM_g$:
$$\MMM_g \into \MMM_g(\lvec{l})\into \wt \MMM_g.$$

\begin{enumerate}
\item[I.] The stack $\MMM_g(\lvec{l})$ is  tame, proper
(separated), smooth, irreducible and of Deligne--Mumford type.
The morphism $\MMM_g(\lvec{l})\to \MMMbar_g$ is finite, flat,
and is an isomorphism on the open dense substack $\MMM_g$.
\item[II.] Any proper substack $\sta X$ of $\widetilde
\MMM_g$ fitting in $\MMM_g \into \sta X\into \wt \MMM_g$ is
isomorphic to $\MMM_g(\lvec{l})$ for a suitable multiindex
$\lvec{l}$.
\end{enumerate}
}

\subsection{The torsor of $r$th roots of a bundle
}
For any line bundle $\sta F$ on the smooth universal curve on $\MMM_g$
whose relative degree is a multiple of $r$
the category
$\MMM_g^{\sta F, r}$
of $r$th roots of $\sta F$ on curves $\sta C\to X$
forms a stack, \'etale and finite
on $\MMM_g$, and equipped with a
torsor structure
under the group stack $\MMM_g^{\Ocal, r}$.

It is well known that $\sta F$ can be written as
a power $\omega^{\otimes k}$ of the relative dualizing sheaf
on the universal curve modulo pullbacks from $\MMM_g$
(Enriques and Franchetta's conjecture \cite{Ha} \cite{Me} \cite{AC}).
Therefore, in view of an extension of
$\MMM_g^{\sta F, r}$ over $\wt{\MMM}_g$,
we focus on the case
$\sta F=\omega^{\otimes k}$ and
we assume
$(2g-2)k\in r\ZZ.$
We compactify $\MMM_g^{\sta F, r}$
in two steps:

\smallskip

\noindent (1) over $\wt \MMM_g$, we construct the stack
parametrizing $r$th roots of
$\sta F=\omega^{\otimes k}$;

\smallskip

\noindent (2) we  restrict such stack
to the compactifications
$\MMM_g(\lvec{l})\subset \wt\MMM_g$ for suitable indexes $\lvec{l}$.

\medskip

For (1), we define the stack $\wt \MMM_g^{\sta F, r}$
of $r$th roots of $\sta F=\omega^{\otimes k}$ on twisted curves
$$\wt \MMM_g^{\sta F, r}=
\{(\stC,\stL,\sta f)\mid \sta f\colon \stL^{\otimes r} \isoto \sta F_\stC\}
_{/\cong}\longrightarrow \wt\MMM_g.$$
Note that $\wt \MMM_g^{\sta F, r}$  can be regarded as the fibred product
$(\wt{\sta {LB}}_g)\ _{\sta k_r}\!\times_{\sta F} \wt \MMM_g$,
where $\wt {\sta {LB}}_g$ is the stack of line bundles on
genus-$g$ twisted curves, $\sta k_r$ is induced by
the $r$th power in $\BBB\GG_m\to \BBB\GG_m$,
and $\sta F$
is regarded as a section $\wt \MMM_g\to \wt {\sta {LB}}_g$.
We show that $\wt \MMM_g^{\sta F, r}$ is a Deligne--Mumford stack,
\'etale on $\wt \MMM_g$.

For (2), we choose a multiindex $\lvec{l}=(l_0,l_1,\dots,l_{\lfloor g/2\rfloor})$
of invertible integers and we consider the restriction
$\MMM^{\sta F,r}_g(\lvec{l})\to \MMM_g(l)$ of $\wt \MMM_g^{\sta F, r}\to \wt \MMM_g$.
In this way, for each $\lvec{l}$, we obtain a stack $\MMM^{\sta F, r}_g(\lvec{l})$ of
$r$th roots of $\sta F$ on $\lvec{l}$-stable curves fibred
over $\MMM_g(\lvec{l})$.
The properness and
the torsor structure are lost
for general choices of $\lvec{l}$ as we already pointed out
in the case
$\lvec{l}=(1,\dots,1)$, which corresponds to stable curves.
The following theorem determines the suitable choices of
$\lvec{l}$.

\vspace{0.2cm}
\noindent{\bf 
\ref{thm:COND} Theorem.}
{\it
For any $\sta F=\omega^{\otimes k}$,
the category  $\MMM^{\sta F, r}_g(\lvec{l})$ is a smooth
Deligne--Mumford algebraic stack, \'etale on $\MMM_g(\lvec{l})$.
\begin{enumerate}
\item[I.] For $\sta F=\Ocal$,
the stack $\MMM_g^{\Ocal , r}(\lvec{l})$ is a finite group stack
if and only if $r$ divides $l_0$.
\item[II.] For $\sta F=\omega$ and $2g-2\in r\ZZ$,
the stack $\MMM_g^{\Ocal , r}(\lvec{l})$ is a finite
group stack and $\MMM_g^{\omega , r}(\lvec{l})$ is a finite  torsor
under $\MMM_g^{\Ocal , r}(\lvec{l})$ if and only if $r$ divides
$$(2i-1)l_i \text { for all $i$}.$$
In this way, we obtain several compactifications
of the stack $\MMM_g^{\omega, r}$ of smooth $r$-spin curves:
for each $\lvec{l}$ satisfying $l_i(2i-1)\in r\ZZ$,
$$\MMM_g^{\omega, r}(\lvec{l})\to \MMM_g(\lvec{l})$$
is the finite torsor of \emph{$r$-spin $\lvec{l}$-stable curves}.
\item[III.] More generally, for $\sta F=\omega^{\otimes k}$  and $(2g-2)k\in r\ZZ$,
the stack $\MMM_g^{\Ocal , r}(\lvec{l})$ is a finite group stack
and $\MMM_g^{\sta F , r}(\lvec{l})$ is a finite  torsor
under $\MMM_g^{\Ocal , r}(\lvec{l})$ if and only if
$r$ divides
$$l_0 \quad \text{and} \quad (2i-1)kl_i, \text { for $i>0$}.$$
\end{enumerate}
}

\medskip

The fact that these compactifications allow a natural
extension of the torsor  structure of $r$th roots defined on the initial
uncompactified moduli stack is an improvement in
its own right. We further mention some concrete situations
 in enumerative geometry where this construction is useful:

\smallskip

\noindent \textbf{Gromov--Witten theory.}
In actual calculations of
enumerative geometry of curves, the main advantage
of our description of $r$th roots via this new notion of stability is
the generalization to $r$th roots of the
classical tools employed
for stable curves.
As an example, in \cite{Ch_Tow}, we illustrate how the
Grothendieck Riemann--Roch formula
allows concrete calculations of
the genus-$g$ Gromov--Witten invariants of the stack $[\CC^2/G]$, where
$G$ is a cyclic subgroup of $SL_2(\CC)$.
This calculation is the subject of the crepant resolution
conjecture, see \cite{BG} for a statement
and \cite{CCIT} for recent progress in genus $0$.
\smallskip

\noindent \textbf{Tensor products of $r$th roots.}
In \cite[Rem.~4.11]{JKV2}
the authors point out that the natural isomorphism
$\MMM^{\sta F,r_1r_2}_{g}\cong \MMM^{\sta F,r_1}_{g}\times_{\MMM_{g}} \MMM^{\sta F,r_2}_{g},$
for relatively prime indexes $r_1$ and $r_2$,
does not extend to the boundary.
In Proposition \ref{pro:r1r2}
we show that the new compactification allows us to
extend the equivalence $L\mapsto (L^{\otimes r_2},L^{\otimes r_1})$
to the boundary over the category of twisted curves with
stabilizers of order $r_1r_2$ on all nodes ($r_1r_2$-stable curves).

In fact, in this way, we get smooth
compactifications of the moduli functor of
$m$-tuples of spin structures of orders $r_1, \dots, r_m$.
These moduli stacks are
used in \cite[Thm.~6.2]{JKV2} in order
to describe tensor products of Frobenius manifolds
(note that, because of ramifications, taking the fibred product
of several Jarvis's compactifications
as in \cite[\S4.2]{JKV2} does not yield a smooth
compactification).

\smallskip

\noindent \textbf{Counting boundary points.}
The description of the geometric points of the boundary locus becomes
straightforward: for the indexes $\lvec{l}$ defined in the previous theorem,
the boundary points are simply
represented by $\lvec{l}$-stable curves with their $r^{2g}$ distinct $r$th roots.
We illustrate in Example \ref{exa:roots} that this improves our
understanding of the enumerative geometry of $r$th roots:
we show how to count the number of $r$-spin structures on a twisted curve $\stC$ up to
automorphisms of $\stC$.
This leads to a counterexample of Conjecture
4.2.1 of \cite{Ja_Pic}, which states that the Picard group of the
moduli stack of smooth genus-1 $r$-spin structures is finite,
see Example \ref{exa:counter}.

Finally, we point out that the results obtained with
previous compactifications extend easily to  the above stacks $\MMM_g^{\sta F,r}(\lvec{l})$.
In Proposition \ref{pro:compAJ}, we show that there is a surjective
morphism from our compactification
to the preexisting compactification due to Abramovich and Jarvis
and illustrate where this morphism is
not invertible.
In Proposition \ref{pro:wtcc}, we prove that the functor of
\cite{PV} and \cite{Ch} defining the
Witten top Chern class
yields a class in the rational cohomology of the new compactification and
we show that such a class is compatible with previous constructions.

\subsection{Structure of the paper}
In Section \ref{sect:terminology}, we fix our terminology and
prove some preliminary results.

In Section \ref{sect:funct} we prove the main technical results.
The subsection 3.1 is a brief subsection where we
prove that the functor of
$r$th roots of a line bundle on a twisted curve on
a base scheme $X$ is a Deligne--Mumford stack \'etale on $X$.
The subsection 3.2 focuses on the geometric fibres of this functor on $X$:
we work out the Kummer theory of a twisted curve and
compare the long exact Kummer sequence
of a twisted curve to that of its coarse space, Theorem \ref{thm:roots}.
In particular, we state a criterion
for a line bundle $\sta F$ to have $r^{2g}$ $r$th roots on a twisted curve,
see Theorem \ref{thm:rootsnum} and see Figure 3.2 at \ref{point}.

By applying these results, in Section \ref{sect:olsson},
we prove Theorem \ref{thm:olsson}
and Theorem \ref{thm:COND} stated above.
We illustrate these theorems in Example \ref{exa:roots}.
We describe the relation with the previous compactification
of Abramovich and Jarvis, Proposition \ref{pro:compAJ}.

In the Appendix we show that line bundles over twisted curves form a stack.
A more general treatment which extends to coherent sheaves and proves that such  a stack is
algebraic can be found in \cite{Li}.

\subsection{Acknowledgements}
I would like to thank
Arnaud Beauville, Alessio Corti,
Carlos Simpson, Angelo Vistoli, and
Charles Walter for
comments and help, and especially Andr\'e Hirschowitz for
numerous stimulating discussions
and for careful reading of preliminary versions of this paper.

\section{Terminology and preliminaries}\label{sect:terminology}
\subsection{Schemes}
We fix an integer $r > 0$,
and throughout this paper we will consider only schemes over \linebreak $\Spec \ZZ[1/r]$.

\subsection{Stacks}
\noindent\textbf{Terminology and generalities.}
Our general reference  is \cite{LM}.
An \emph{algebraic stack} is a stack satisfying Artin's definition \cite{Ar}.
Stacks in the sense of Deligne and Mumford \cite{DM}
will be called \emph{Deligne--Mumford stacks}.
When working with algebraic stacks with finite diagonal,
we use Keel and Mori's Theorem
\cite{KM}: there exists an algebraic space $\coa{\sta X}$ associated to
the stack $\sta X$
and a morphism $\pi_{\sta X}\colon \sta X\to \coa{\sta X}$
(or simply $\pi$)
which is universal with respect
to morphisms from $\sta X$ to algebraic spaces. We
refer to $\coa{\sta X}$ as the \emph{coarse space}.
In this way we have a functor (and in fact a $2$-functor)
associating to any morphism between this type of stacks
$\sta f\colon \sta X\to \sta Y$ the unique morphism
between the corresponding coarse algebraic spaces
$\coa{\sta f}\colon \coa{\sta X}\to \coa{\sta Y}$
satisfying $\sta f\circ \pi_{\sta Y}=
\pi_{\sta X}\circ \coa{\sta f}$.

We refer to \cite{Br} for the notion of
\emph{group stack} $\sta G\to \sta X$.
We say that there is an action of the group stack $\sta G\to \sta X$
with product $\sta m_{\sta G}$ and unit object $\sta e$
on $\sta T\to \sta X$
if there is a morphism of stacks
$\sta m\colon \sta G\times _{\sta X}\sta T\to \sta T$
and homotopies
$\sta m\circ (\sta m_{\sta G}\times \id_{\sta T})
\Rightarrow \sta m(\id_{\sta G}\times \sta m)$
and $\sta m\circ(\sta e\times \id_{\sta T})\Rightarrow \id_{\sta T}
$
satisfying the associativity constraint \cite[6.1.3]{Br} and
the compatibility constraint \cite[6.1.4]{Br}.
The morphism $\sta T\to\sta X$ is a \emph{torsor}
if the morphism
$$\sta m\times \sta {pr}_2\colon \sta G\times _{\sta X}\sta T\to \sta T\times _{\sta X} \sta T$$
is an isomorphism of stacks and $\sta T\to \sta X$
is flat and surjective.

\vspace{0.3cm}

\noindent\textbf{Morphisms of stacks.}
We often need to consider $2$-categories in which the
objects are algebraic stacks, the functors between two stacks are regarded
as $1$-morphisms, and the natural transformations are regarded as
$2$-morphisms.

The situation is often simplified
by the following criterion
showing that certain morphisms between stacks
have only trivial $2$-automorphisms. In particular
this criterion applies to
morphisms between twisted curves.
\begin{lem}[Abramovich and Vistoli, {\cite[Lem.~4.2.3]{AV}}]\label{lem:AV}
Let $\sta f \colon \sta X \to \sta Y$
be a representable morphism of Deligne--Mumford stacks over a scheme $S$.
Assume that there exists a dense open representable substack
(i.e. an algebraic space) $U \subset  \sta X$ and an open
representable substack $V\subset Y$ such that $\sta f$ maps $U$ into $V$.
Further assume that the diagonal $\sta Y \to \sta Y \times_S \sta Y$
is separated. Then any automorphism of $\sta f$ is trivial.
\qed
\end{lem}


\noindent\textbf{Stabilizer of a geometric point of a stack.}
Let $\sta X$ be an algebraic stack.
A geometric point $\sta p\in \sta X$
is an object  $\Spec k\to \sta X$, where $k$ is algebraically closed.
We denote by $\Aut(\sta p)$ the automorphism group
of $\sta p$ as an object of the fibred category $\sta X_\sta p$.
We refer to $\Aut(\sta p)$ as
the \emph{stabilizer} of $\sta p$.

\vspace{0.3cm}

\noindent\textbf{Local pictures.}
We often need to describe stacks and morphisms between stacks locally
in terms of explicit equations.
We adopt the following standard convention, which avoids repeated mention
of strict henselization \cite[\S1.5]{ACV}.

Let ${\sta X}$ and $\sta U$ be algebraic stacks and
let $\sta x\in \sta X$ and $\sta u\in \sta U$ be geometric points.
We say ``the local picture of ${\sta X}$
at $\sta x$ is
given by ${\sta U}$ (at $\sta u$)'' if
there is an isomorphism between the
strict henselization ${\sta X}^{\text{sh}}$
of ${\sta X}$ at $\sta x$ and the strict henselization
${\sta U}^{\sh}$
of ${\sta U}$ at $\sta u$.

If $\sta f\colon \sta X\to \sta Y$ and
$\sta g\colon \sta U\to \sta V$ are morphisms of stacks
and $\sta x$ and $\sta u $ are geometric points
in $\sta X$ and $\sta U$, we say ``the local picture of
$\sta X \to \sta Y$ at $\sta x$ is given by $\sta U\to \sta V$ (at $\sta u$)''
if there is an isomorphism between the strict henselization
${\sta f}^{\text{sh}}\colon \sta X^{\sh}\to \sta Y^{\sh}$
of $\sta f$ at $\sta x$
and the strict henselization
${\sta g}^{\text{sh}}\colon \sta U^{\sh}\to \sta V^{\sh}$
of $\sta g$ at $\sta u$. This convention
allows local descriptions of diagrams of morphisms
between stacks; in particular it allows local description
of  group actions $G\times \sta X\to \sta X$ and of $G$-equivariant
morphisms.

\vspace{0.3cm}

\noindent\textbf{The construction of Cadman, Matsuki, Olsson, and Vistoli.}
\label{sect:MOC}
For any smooth scheme $X$ and smooth effective  divisor $D$ in $X$
we present a stack-theoretic modification of $X$, which contains
$X\setminus D$ as a dense open representable substack
and has coarse space $X$.

In \cite{Ca}, Cadman provides the following construction.
\begin{defn}
Let $X$ be
a smooth scheme $X$, let
$D$ be a smooth effective Cartier divisor in $X$,
and let $l$
be a positive integer invertible on $X$.
The category $X[D/l]$
is formed by objects
$(S, M, j, s)$, where
\begin{enumerate}
\item $S$ is an $X$-scheme $S\to X$;
\item $M$ is a line bundle on $S$;
\item $j$ is an isomorphism between $M^{\otimes l}$ and the pullback of $\Ocal (D)$ on $S$;
\item $s$ is a section $s\in \Gamma(S,M)$ such that
$j(s^{\otimes l})$ equals the tautological section
of $\Ocal(D)$ vanishing along $D$.
\end{enumerate}
The morphisms are defined in the obvious way.
\end{defn}
This definition
yields a Deligne--Mumford stack
$X[D/l]$, with coarse space $X$.
The morphism $\pi \colon X[D/l]\to X$ is an isomorphism over $X\setminus D$: we have
$$X\setminus D\into X[D/l]\xrightarrow{\ \pi\ } X,$$
where $X\setminus D$ is dense in $X[D/l]$ and
$\pi$ is finite and flat.
Note that $X[D/l]$ is
equipped with a tautological line bundle
$\sta M$ and an isomorphism
\begin{equation}
\label{eq:taut}\sta M^{\otimes l}\isoto \pi^*\Ocal_X(D).\end{equation}

A special case of this construction was first introduced
by Abramovich, Graber, and Vistoli \cite[3.5.3]{AGV}
(the idea is attributed to Vistoli, see \cite[3.5]{AbrNotes}).
In the existing literature, two different and
compatible definitions  can be found: see
Matsuki and Olsson \cite{MO} and Cadman \cite{Ca} for the above definition
(see \cite[2.4.5]{Ca} for the compatibility
between the two constructions).

It is natural to try and generalize this construction.
If $D_1, \dots, D_n$ are distinct
smooth effective divisors with normal crossings,
for any positive integers $l_1,\dots, l_n$ invertible on $X$,
we write
\begin{equation}\label{eq:conv}
X[D_1/l_1+\dots +D_n/l_n]:= X[D_1/l_1]\times_X \dots \times_X X[D_n/l_n].
\end{equation}
In this way, $X[D_1/l_1+\dots +D_n/l_n]$ is a smooth Deligne--Mumford stack.

\begin{gener}\label{gen:log}
The definition provided in \cite{MO} by Matsuki and Olsson
generalizes the above definition of $X[D_1/l_1+\dots +D_n/l_n]$
to the case where singular divisors $D_i$ occur, and
the local picture of $D_i$ at each point is the union of
smooth divisors with normal crossings (a normal crossings divisor).
This extension involves the notion of logarithmic structures
in the sense of Fontaine and Illusie or the use of \'etale descent.
The output is again a smooth stack, see \cite[Thm.~4.1]{MO}.
\end{gener}

\begin{exa}\label{affineMOC}
Consider the
affine space $X=\Spec R$ and
the divisor $D=\{t=0\}$ for $t\in R$.
Assume that $l$ is invertible in $R$.
The stack $X[D/l]$ is the quotient stack
$[\Spec \wt R/\pmmu_l]$, where
$\wt R=\Spec R[\wt t]/(\wt t^l-t)$ and
$\pmmu_l$ acts on $\wt R$ as $g\cdot \wt t=g^{-1}\wt t$ and
fixes $R$.
\end{exa}

\begin{exa}\label{exa:onedim}
Let $\sta X$ be a Deligne-Mumford stack,
whose coarse space is
a proper, regular, and reduced curve $\coa{\sta X}$ over a field $k$
and whose geometric points have trivial stabilizers except for
a finite number of distinct points
$\sta p_1,\dots ,\sta p_n$ with stabilizers of
order $l_1,\dots, l_n$ (which we assume invertible).
In that case, we have (see \cite[Exa.~3.7]{Ca})
$$\textstyle {\sta X\cong \coa{\sta X}[p_1/l_1+\dots +p_n/l_n]}$$
where $p_i=\coa {\sta p_i}\in \coa{\sta X}$
is the point corresponding to $\sta p_i\in \sta X$.
The local picture at $\sta p_i\colon \Spec k \to \sta X$ is given by $S[0/l_i]$, where
$S=\Spec k[z]$ and $0$ denotes $\{z=0\}$.
In this way there is an isomorphism
between
$\Aut(\sta p_i)$ and
cyclic group of $l_i$th roots of unity of $\GG_m$.

\begin{rem}\label{rem:cangener}
At a point $\sta p_i\colon \Spec k\to \sta X=\coa{\sta X}[p_1/l_1+\dots + p_n/l_n]$
the projection
$\sta {pr}_2$ from $\Spec k \,_{\coa{\sta p_i}}\!\times_\pi \sta X$ to $\sta X$
induces a canonical embedding $\sta j_i\colon \sta B(\Aut(\sta p_{i}))\into \sta X$.
The group $\Pic(\sta B(\Aut \sta p_{i}))$
is a cyclic group of order $l_i$.
The tangent space
$\sta T_i$ at $\sta p_i$ is a representation of
$\Aut(\sta p_i)$ and is a canonical generator of
$\Pic(\Aut \sta p_{i})$. In this way
$\Pic(\Aut \sta p_{i})$ is canonically isomorphic to $\ZZ/l_i\ZZ$.
\end{rem}
\end{exa}

\begin{rem}
On a Deligne--Mumford stack of dimension $1$, the
degree of a line bundle $\sta F$ is the degree of
the first Chern class $c_{1}(\sta F)$
in the rational Chow ring $\in A^1(\sta X)_\QQ$.
Therefore $\deg(\sta F)$ is a rational number,
see \cite[3.3]{Kr} and references therein.
In the case of $C[p/d]$,
we have the following proposition.
\end{rem}
\begin{pro}\label{pro:deg}
Let $C$ be a proper, regular, and reduced curve over a field $k$.
For a closed point $p\in C$ and $l$ a positive integer,
consider $\pi \colon C[p/l]\to C$.
We write $\Gamma$ for the stabilizer of the point of $C[p/l]$ over $p$ and
write $\sta j\colon \sta B\Gamma\to C[p/d]$
for the canonical embedding of Remark \ref{rem:cangener}.
For any line bundle $\sta F$ on $C[p/d]$, the degree of
$\sta F$ belongs to $\frac{1}{l}\ZZ$ and we have
$$l\deg(\sta F)\equiv \sta j^*(\sta F) \in \ZZ/l\ZZ,$$
where we used the canonical identification of $\Pic (\sta B\Gamma)$
with $\ZZ/l\ZZ$ of Remark \ref{rem:cangener}.
\end{pro}
\begin{proof}
Let $h\in \{0,\dots , l-1\}$ be the integer satisfying
$\sta j^*(\sta F)=\sta T^{\otimes h}$, for $\sta T$
the tangent space at $\sta p$.
The tautological line bundle $\sta M$ on
$C[p/d]$ satisfies $\sta M^{\otimes l}\cong \pi^*\Ocal_C(p)$,
\eqref{eq:taut}.
Note that $\sta j^*\sta M=\sta T$.
Hence, $\sta F\otimes (\sta M^{\otimes h})^\vee=\pi^*\Ocal_C(D)$
for an integral divisor $D$ on $C$, and we have
\begin{align*}
l\deg(\sta F)&=
l\deg(\sta M^{\otimes h}\otimes \pi^*\Ocal_C(D))\\
&=\deg(\sta M^{\otimes hl}\otimes \pi^*\Ocal_C(lD))\\
&=\deg(\pi^*\Ocal(h[p])\otimes \pi^*\Ocal_C(lD))\\
&=\deg_C(\Ocal(h[p])\otimes \Ocal_C(lD))\\
&=h+l\deg D\\
&\equiv h \mod l\\
&\equiv \sta j^*(\sta F) \mod l.
\end{align*}
\end{proof}
\subsection{Stable curves}
We recall some standard definition for sake of clarity and
for preparing the definition of twisted curve in \S2.4.
For any integer  $g\ge 2$, a \emph{stable} curve
of genus $g$ on a scheme $X$ is
a proper and flat morphism
$C\to X$ satisfying the following conditions:
\begin{enumerate}
\item each geometric fibre $C_x$ over $x$ in $X$ is reduced, connected,
has dimension one, and has only ordinary double points (which we
usually call \emph{nodes});
\item the dualizing sheaf of $C_x$  is ample;
\item $\dim_{k(x)}H^1(C_x,\Ocal_{C_x})=g$.
\end{enumerate}

A stable curve $C$ over an algebraically closed field $k$ yields a natural
combinatorial object: the \emph{dual graph}.
The dual graph $\Lambda$
of $C$
is the graph whose
set of vertices
$V$ is the set of irreducible components
of $C$ and whose set of edges
$E$ is the set of nodes of $C$. An edge
connects the vertices corresponding
to the irreducible components containing the
two branches of the nodes.
If an orientation of $\Lambda$ is fixed, we have a chain complex
$\Ccal_\bullet(\Lambda,\ZZ/r\ZZ)$ with differential
\begin{align*}
\partial \colon (\ZZ/r\ZZ)^E&\xrightarrow{\ } (\ZZ/r\ZZ)^V,
\end{align*}
where the edge starting at $v_-$ and ending at $v_+$
is sent to the $0$-chain $[v_+]-[v_-]$.
We say that a node $e$ of a stable curve $C$
is \emph{separating} if by normalizing $C$ at
the point $e$ we obtain two disjoint components.

If we assign to each vertex $v\in V$
the genus $g_v$ of the connected
component of the normalization of $C$ corresponding to the irreducible
component attached to $v$, the (arithmetic) genus of $C$ can be read off
from $\Lambda$
and the function $v\mapsto g_v$. Indeed, we have
\begin{equation}\label{genusgraph}
g(C)=b_1+\textstyle{\sum_v g_v}
\end{equation}
where $b_1=1-\size{V}+\size{E}$ is the first Betti number of $\Lambda$.

\subsection{Twisted curves}\label{notn:twisted}
We recall the notion of twisted curve due to Abramovich and Vistoli,
see \cite{AV} or the remarks
below for some slight generalizations.
\begin{defn}\label{defn:twisted}
A \emph{twisted curve} of genus $g$
on a scheme $X$ is a proper and flat
morphism of tame stacks $\stC\to X$, for which
\begin{enumerate}
\item
the fibres are purely $1$-dimensional
with at most nodal singularities,
\item the coarse space is a stable curve
$\coa{\sta C}\to X$ of genus $g$;
\item the smooth locus
$\stC^{\rm sm}$ is an algebraic space;
\item
the local picture
at a node
is given by $[U/\pmmu_l]\to T$, where
\begin{enumerate}
\item[$\bullet$] $T=\Spec A$,
\item[$\bullet$] $U=\Spec A[z,w]/(zw-t)$ for some $t\in A$, and
\item[$\bullet$] the action of $\pmmu_l$ is given by
$(z,w)\mapsto (\xi_lz,\xi_l^{-1}w)$.
\end{enumerate}
\end{enumerate}
(Recall that the tameness condition on $\stC$ means that
for every geometric point
$\sta p\colon \Spec k \to \sta C$
the group $\Aut(\sta p)$
has order prime to the characteristic of the
algebraically closed field $k$.)
\end{defn}

\begin{rem}[unbalanced twisted curves]
In the existing literature
twisted curves satisfying the above
local condition (4)
are often called \emph{balanced}, \cite{AV}.
We drop the adjective balanced, because
we never consider unbalanced twisted curves.
\end{rem}

\begin{rem}[twisted curves with smooth stack-theoretic points]
A further natural generalization, which is usually
given in the definition of a twisted curve
(see \cite[Def.~4.1.2]{AV}),
consists of allowing
nontrivial stabilizers on smooth points.
By \cite[Thm~4.1]{Ca} and \cite[Thm~1.8]{Ol},
giving
a stack-theoretic curve with a smooth point whose stabilizer
has order $l$ is equivalent
to assigning an invertible integer $l$,
a twisted curve $\stC\to X$ in the sense of
Definition \ref{defn:twisted}, and
a section $\sigma\colon X\to \stC$ in the
smooth locus.
(With the notation introduced above, the equivalence is essentially
$(\stC, \sigma,l)\mapsto {\stC}[\sigma(X)/l]$.)
\end{rem}
\begin{rem}[unstable coarse space]
Twisted curves are generally defined without
imposing the
condition of stability on the coarse space.
In Olsson's paper \cite{Ol}, it is shown that
these stack-theoretic curves
form an algebraic stack in the sense of Artin's definition.
If we require that the coarse space
$\coa{\sta C}$ is stable, we get
a Deligne--Mumford stack. See Theorem \ref{thm:Olsson}.
\end{rem}
A $1$-morphism $\stC\to X$ to $\stC'\to X'$
between twisted curves is a morphism
of stacks $\sta m$ fitting in the fibre diagram
\[\xymatrix@R=0.1cm{\stC\ar[rr]^{\sta m}\ar[dd]&&\stC'\ar[dd]\\
&\square&\\
X\ar[rr]&&X'. }\]
By Lemma \ref{lem:AV}, there is at most one
natural transformation,
which identifies two morphisms. Therefore,
we obtain a \emph{category} by considering
morphisms as
$1$-morphisms up to
base-preserving natural transformations.

In \cite[Thm.~1.9]{Ol}, Olsson proves that such a category is
an algebraic stack. In particular, see \cite[Rem.~1.10]{Ol}, he shows
how the versal deformation space of a twisted curve
$\stC$ over a field $k$ relates
to the versal deformation space
$\Spec I$ of the coarse space
$\coa{\sta C}$. Let $\sta e_1,\dots ,\sta e_m$ be the nodes of $\stC$
and $\coa{\sta e_1}, \dots, \coa{\sta e_m}$
the corresponding nodes of $\coa{\sta C}$.
Let $l_i$ be the order of the automorphism group
of $\sta e_i$. Let \begin{equation}
\label{eq:paramt}
D_i=\{t_i=0\}\end{equation} be the divisor
classifying deformations of $\coa{\sta C}$
on which the node $\coa{\sta  e_i}$ persists:
\begin{theor}[Olsson, {\cite[Thm.~1.9]{Ol}}]\label{thm:Olsson}
The category of twisted curves
$\wt \MMM_g$ in the sense of Definition
\ref{defn:twisted}
is a smooth Deligne--Mumford stack
and the versal deformation space of
a twisted curve $\stC\to \Spec k$ is given by
\begin{equation}
\label{eq:versal}
\wt I:=I[z_1, \dots, z_m]/(z_1^{l_1}
 - t_1, \dots , z_m^{l_m} - t_m),
\end{equation}
where $t_1,\dots, t_m \in I$ are the parameters satisfying \eqref{eq:paramt}
and $l_1, \dots, l_m$ are the orders of the stabilizers of the corresponding nodes.
\qed\end{theor}
The group of automorphisms of a twisted curve
acts naturally      on the versal deformation space
of a twisted curve.
Using \cite{ACV}, we describe the automorphisms
of a twisted curve. Then, in Remark \ref{rem:action},
we illustrate the action on the versal deformation.

\vspace{0.3cm}

\noindent\textbf{Automorphisms of twisted curves.}
Let $\pi\colon \stC\to \coa{\sta C}$ be a twisted curve over
an algebraically closed field $k$.
In \cite[Prop.~7.1.1]{ACV}
the group
$\Aut(\stC,\coa{\sta C})$ of automorphisms
of $\stC$ that fix the coarse space $\coa{\sta C}$
is explicitly calculated:
\begin{theor}[Abramovich, Corti, Vistoli, {\cite[Prop.~7.1.1]{ACV}}]
\label{thm:aut_tw}
For a twisted curve over an algebraically closed field, we have
an isomorphism $$\Aut(\stC,\coa{\sta C})\cong\pmmu_{l_1}\times
\dots\times
\pmmu_{l_m},$$
where ${l_1},\dots,{ l_m}$ are the orders of the
stabilizers at the nodes $\sta e_1,\dots,\sta e_m$.
\qed\end{theor}
\begin{rem}\label{rem:lpghost}
In fact \cite[Thm.~7.1.1]{ACV}
shows that we can choose
generators $\sta g_1,\dots, \sta g_m$
of $\Aut(\stC,\coa{\sta C})$ such that the restriction of
$\sta g_i$ to ${\stC\setminus \{\sta e_i\}}$ is the identity,
and the local picture at $\sta e_i$ is given by
\begin{align*}
k[z,w]/(zw)&\to k[z,w]/(zw) \\
(z,w) &\mapsto (z,\xi_{l_i}w),
\end{align*}
where $\xi_{l_i}$ is a primitive ${l_i}$th root of unity.
(Local pictures are given up to natural transformations, and the
$1$-automorphism above is in fact, locally
on the strict henselization, $2$-isomorphic to
$(z,w)\mapsto (\xi_{l_i}^{a}z,\xi_{l_i}^{b}y)$
for any $a, b\in\{0,\dots, l_i-1\}$ satisfying $a+b=1\mod l_i$.)
\end{rem}

\begin{rem}\label{rem:action}
The group $\Aut(\stC,\coa{\sta C})$,
with the basis given above, acts on the versal deformation
space \eqref{eq:versal}
as $(\sta g_{1},\dots, \sta g_{m})z_i=\xi_{l_i}z_i$ (see
\cite[Lem.~5.3]{Ol}).
\end{rem}

\subsection{Line bundles on twisted curves}
In view of the study of $r$th roots in the Picard group
we review some facts on line bundles on twisted curves.
\begin{pro}\label{pro:omega}
For any twisted curve $\sta f\colon \sta C\to X$ with coarse space
$\coa{\sta f}\colon \coa{\sta C}\to X$ the pullback via
$\pi\colon \sta C\to \coa{\sta C}$ of the
relative dualizing sheaf
of ${\coa{\sta f}}$ is the relative dualizing sheaf
of ${\sta f}$.
\end{pro}
\begin{proof}
The natural homomorphism $\pi^*\omega_{\abs{\sta f}}\to \omega_{\sta f}$
is an isomorphism on the smooth locus.
At the nodes the local picture of $\sta C\to \coa{\stC}$
is
\begin{align*}
U=\Spec A[z,w]/(zw-t)&\to V=\Spec A[x,y]/(xy-t^l)\\
(z,w)&\mapsto (z^l,w^l)
\end{align*}
with the stabilizer of the node $\pmmu_l$ acting as
$\xi_l\cdot (z,w)=(\xi_lz,\xi_l^{-1}w)$.
The local generator of the relative dualizing sheaf of $\coa{\sta f}$ is
$$\frac{dz^l}{z^l}-\frac{dw^l}{w^l}$$
and is sent to the local generator $dz/z-dw/w$ of $\omega_{\sta f}$ (recall
that $d(z^l)=lz^{l-1}dz$ and $d(w^l)=lw^{l-1}dw$).
\end{proof}

We now fix a twisted curve over an algebraically closed field
$k$ and describe the action by pullback
of $\Aut(\stC)$ on $\Pic(\stC)$.
We take one of the automorphisms belonging to $
\Aut(\sta C, \sta C\setminus \{\sta e\})\subseteq \Aut(\sta C, \coa{\stC})$
fixing the entire twisted curve away from a node $\sta e$ with stabilizer of
order $l$. We describe explicitly such automorphism
by choosing for the rest of this section a primitive $l$th root of unity $\xi$.
The local picture of  $\stC$  at $\sta e$
is  given by $[V/\pmmu_l]$, where $V$ equals $\{z_+z_-=0\}$ and $\pmmu_l$ acts as
$(z_+,z_-)\mapsto (\xi z_+,\xi^{-1}z_-)$.
The automorphism
$\sta g\in \Aut(\stC,\coa{\stC})$ is
chosen in such a way that $\sta C\setminus\{\sta e\}$
is fixed and the local picture at $\sta e$ is
\begin{equation}\label{eq:gloc}
{\sta g}^{\rm sh}\colon (z_+, z_-)\mapsto (z_+,\xi z_-).
\end{equation}
We already noted in Theorem \ref{thm:aut_tw} and
in Remark \ref{rem:lpghost} that all automophisms of $\stC$
are equal (up to natural transformation)
to the composite of pullbacks from $\coa{\stC}$ and products of morphisms
defined in the same way as $\sta g$.

We now take a line bundle $\stL$ on $\stC$ and show in
Proposition \ref{pro:aut_lb} that
pulling back
$\stL$ via $\sta g$ is
the same as tensoring $\stL$ by a
line bundle $\sta T_\sta L$ in the torsion subgroup of the Picard group of $\stC$.
We need to set up some standard notation
identifying the line bundle $\sta T_{\sta L}$.

Consider the homomorphism $$\gamma\colon \GG_m\to \Pic \stC$$
sending $\lambda\in \GG_m$ to the line bundle
of regular functions
$f$ on the partial  normalization at
the node $\sta e$ satisfying
$f(\sta p_+)=\lambda f(\sta p_-)$, where $\sta p_+$ and $\sta p_-$ are the points
of the normalization lifting the node $\sta e$.

Consider the pullback of $\stL$ on the partial normalization, and its restriction
on ${\sta B}(\Aut(\sta p_+))$. We already noted in Remark \ref{rem:cangener}
that the Picard group of ${\sta B}(\Aut(\sta p_+))$
is canonically generated by $\sta T_{\sta p_+}$. Therefore,
$\sta L$ determines  on ${\sta B}(\Aut(\sta p_+))$ a power of $\sta T_{\sta p_+}$;
we denote its exponent by
$$\mult_{\sta p_+}{\sta L}\in \{0, ... , l-1\}.$$

Note that $\gamma$ and $\mult_{\sta p_+}\sta L$ are transformed into
 $\la\mapsto \gamma(\la^{-1})$ and into
 $l-\mult_{\sta p_+}\sta L$ (reduced modulo $l$) if
 we interchange the notations $\sta p_+$ and $\sta p_-$. This implies that
the line bundle
$$\sta T_{\sta L}=\gamma(\xi^{\mult_{{\sta p}_+} \sta L})\in \Pic(\sta C)$$
only depends on $\sta L$ and does not depend on the notation $\sta p_+$, $\sta p_-$.

\begin{pro}\label{pro:aut_lb}
For any line bundle $\stL$ on $\stC$ over an algebraically closed  field $k$
and for $\sta g\in \Aut(\sta C, \sta C\setminus\{\sta e\})$ satisfying \eqref{eq:gloc}
we have
\begin{equation}
\label{eq:adjust}
\sta g^*\stL\cong \stL\otimes \sta T_{\sta L}.
\end{equation}\qed
\end{pro}
\begin{proof}
The local picture of $\stL\to \stC$
at the point of the zero section over $\sta e$ is
$$\sta W=[(V\times \Af^1)/\pmmu_{l}]\to [V/\pmmu_l]$$
where $V$ is $\{z_+z_-=0\}$ as above, and
$\xi\in \Aut(\sta e)$ acts on $((z_+,z_-),\la)\in V\times \Af^1$
as
\begin{equation}\label{eq:param_h}
((z_+,z_-), t)\mapsto
((\xi z_+,\xi^{-1}z_-), \xi^{\mult_{{\sta p}_+}\sta L} t ). \end{equation}

By tensoring with a suitable
element of $\pi^*\Pic\coa{\sta C}$,
we can restrict to the case of a
line bundle $\stL$ on $\stC$ which
is trivial on $\stC\setminus\{\sta e\sta\}$.
In this way, for $V^{\times}=V\setminus (0,0)$,
we can regard $\stL$ as the datum
of a line bundle $\sta W$ on $[V/\pmmu_{l}]$ alongside with
an isomorphism $\Phi$ between $\sta W\rest {V^\times}$
and the submodule of
$\Ocal_{V^\times}= k[z_+,z_+^{-1}]\oplus
 k[z_-,z_-^{-1}]$ invariant under the action of $\pmmu_l$.
By \eqref{eq:param_h}, the line bundle $\sta W$ is trivial on $V$ and
$\pmmu_{l}$-linearized
by the character $\xi \mapsto\xi^h$ for $h={\mult_{{\sta p}_+}\sta L}$.
Note that invariant sections of
$\sta W\rest {V^\times}$ form a module
$z_+^{l-h} k[z_+^{l},z_+^{-{l}}]\oplus
z_-^{h} k[z_-^{l},z_-^{-l}]$.
Pulling back via $\sta g:(z_+, z_-)\mapsto (z_+, \xi z_-)$ changes
$\Phi$ by multiplication by
$(1\oplus \xi^h)$ in $k[z_+,z_+^{-1}]\oplus
k[z_-,z_-^{-1}]$.
In this way, we only change the descent datum
along the partial desingularization at $\sta e$
by multiplication by  $\xi^h$
on the branch $(z_+=0)$. This amounts to tensoring as in  \eqref{eq:adjust}.
\end{proof}

On a stable curve $C\to \Spec k$, the total degree of a
line bundle is the sum of the
degrees on the connected components of the normalization.
We extend this definition to twisted curves.
\begin{pro}
The total degree of a line bundle $\sta F$ on a twisted curve $\sta C$
over $k$ is an integer.\end{pro}
\begin{proof}
We only need to prove the claim in the case
when $\sta F$ is trivial on $\stC\setminus\{\sta e\}$
and nontrivial on the node $\sta e$.
Then, in order to calculate the total degree, we can regard $\sta F$ as
a line bundle on a twisted curve
$\stC$  with trivial stabilizers
on every node except $\sta e$. The
normalization of $\stC$ is a smooth stack $\sta C^\nu$
as in Example \ref{exa:onedim}:
$\sta C^\nu=X[p_1/l+p_2/l]$
for $p_1,p_2\in X$ and $l\in \ZZ_{\ge 1}$. We calculate the total degree using Proposition
 \ref{pro:deg}.

First, assume that the partial normalization of $\stC$ at $\sta e$
is the disconnected stack $\stC^\nu=\sta D_1\sqcup \sta D_2$.
By Proposition \ref{pro:deg}, the total degree of $\sta F$ is
in $\frac{1}{l}\ZZ$.
We need to show that the total degrees of the
restrictions $\sta F_1$ and $\sta F_2$ on $\sta D_1$ and $\sta D_2$
satisfy \begin{equation}
\label{eq:congruence_relation}
l \deg(\sta F_1)+l\deg(\sta F_2)\equiv 0\mod l.
\end{equation}
By Proposition \ref{pro:deg}, this amounts to showing that the
pullbacks of $\sta F_i$ with respect to \linebreak
$\sta B\Aut(\sta p_i) \to \sta D_i$
yield inverse characters of
for $i=1$ and $2$.
This follows using the local picture given above, \eqref{eq:param_h}.

Finally, if the partial normalization
$\sta h\colon {\sta C}^\nu \to \stC$ is
connected we have two distinct points
$\sta p_1$ and $\sta p_2$
with nontrivial stabilizer of order $d$ lying over $\sta e$.
We can define two line bundles
$\sta F_1$ and $\sta F_2$ on ${\sta C}^\nu$
such that  $\sta F_i$ is trivial at $\sta p_i$,
and $\sta F_1\otimes \sta F_2=\sta h^*\sta F$.
Then, \eqref{eq:congruence_relation} holds (with the
same proof) and this implies the claim.
\end{proof}

\section{The functor of $r$th roots of a line bundle}\label{sect:funct}
This section is divided in two subsections: \S\ref{sect:funct}.1 and \S\ref{sect:funct}.2.
\begin{enumerate}
\item Relying on the results of \cite{Li} or on the appendix, this first part
recalls briefly that $r$th roots of a line bundle form a Deligne--Mumford stack.
We consider a twisted curve
$\stC\to X$
and the functor ${\sta F^{1/r}}$ of
$r$th roots of $\sta F$, a line bundle on $\stC$
whose relative degree is a multiple of $r$.
We show that it is a $\pmmu_r$-gerbe over a scheme \'etale on $X$.
We notice that the stack of $r$th roots ${\sta F^{1/r}}$
is finite as soon as the geometric fibre on $X$ is constant.
\item
In Section \ref{sect:kummer}, we study the geometric
fibres
of $\sta F^{1/r}\to X$
by calculating the cohomology of
the Kummer sequence
$$1\to \pmmu_r\to \GG_m\xrightarrow{\ r\ } \GG_m\to 1$$
for a twisted curve $\stC$ over an algebraically closed field.
The main result is Theorem \ref{thm:roots},
where the endomorphism $\stL\mapsto \stL^{\otimes r}$
of $\Pic \stC$ is inscribed in a diagram of exact sequences.
In Corollary \ref{cor:root_tw}, we apply the result to
twisted curves (we draw
the diagram for a twisted curve in Figure 3.2 at \ref{point}).
Finally, in Theorem \ref{thm:rootsnum}, we deduce the
numerical criterion classifying
line bundles on twisted curve having exactly $r^{2g}$ roots.
Let us record straight away a
consequence of Theorem \ref{thm:roots} in the context
where $r$ divides $\size{\Aut(\sta e)}$ for every node $\sta e$.
\end{enumerate}
\begin{cor}\label{cor:rdividesl}
Let $\pi\colon \stC\to \coa{\stC}$ be
a twisted curve of genus $g$ over an algebraically closed field.
There is an exact sequence
$$1\to \Pic \coa{\stC}\xrightarrow{\pi^*}
\Pic \stC\to \textstyle{ \prod_{\sta e\in E} }
\Pic \sta B(\Aut (\sta e))\to 1,$$
where $E$ is the set of nodes
$\sta e$ in $\stC$.
Furthermore, as soon as $$\size{\Aut(\sta e)}\in r\ZZ \quad\quad \quad \forall \sta e\in E $$
we can write
$\textstyle{ \prod_{E} }
\Pic \sta B(\Aut (\sta e))_r\cong \Ccal_1(\Lambda,\ZZ/r\ZZ)$
where $\Lambda$ denotes the dual graph of $\stC$
and we have an exact sequence
\begin{equation}\label{eq:gist}
1\to (\Pic \coa{\stC})_r
\xrightarrow{\pi^*}
(\Pic \stC)_r\to \Ccal_1(\Lambda,\ZZ/r\ZZ)\xrightarrow {\partial}
\Ccal_0(\Lambda,\ZZ/r\ZZ)\xrightarrow {\varepsilon}\ZZ/r\ZZ\to 1
\end{equation}
where $(\Pic \coa{\stC})_r$ and $(\Pic \stC)_r$ denote
the $r$-torsion subgroups of the Picard groups, $\partial$ is
the boundary homomorphism with respect to a chosen orientation of
$\Lambda$, and $\varepsilon$ denotes the augmentation homomorphism
sending $(h_v)_V$ to $\sum_V h_v\in \ZZ/r\ZZ$.

As a consequence, for any line bundle
on $\coa{\stC}$, whose total degree is a multiple of $r$,
the pullback on $\stC$ has $r^{2g}$ $r$th roots.
(For $(2g-2)k\in r\ZZ$, this applies in particular to
$\sta F=\omega^{\otimes k}_{\stC}$
by Proposition \ref{pro:omega}.) \qed
\end{cor}


\subsection{The stack of $r$th roots of a line bundle}\label{sect:F/r}
Consider the twisted curve $\sta f\colon \stC\to X$.
We write
$\sta{LB}_{\sta f}$
for the category of line bundles on
base changes $\sta C_S=\sta C\times_X S$
for every $X$-scheme $S$.
More precisely, the \emph{objects} are pairs $(S,\sta M)$, where
$S$ is an $X$-scheme
and $\sta M$ is a line bundle on
$\sta C_S=\sta C\times_X S$. The
\emph{morphisms} $(S,\sta M)\to (S',\sta M')$
are pairs $(m,\sta a)$, where $m\in \mathsf{Hom}_X(S,S')$
and $\sta a$ is an isomorphism of line bundles
$\sta a\colon \sta M\xrightarrow{\sim} \sta M'\times _{S'} S$
on $\sta C_S$.
This fibred category is an algebraic stack and indeed a substack of
the fibred category of coherent sheaves, \cite[Thm.~2.1.1, Lem.~2.3.1]{Li}.
The properties of $\sta {LB}_{\sta f}$ needed here are independently proven in
Appendix A (Proposition \ref{thm:lb} and Remark \ref{rem:algebraic}).

\begin{rem}\label{rem:globalLB}
Note that, in this way, we can
define the algebraic stack $\wt{\sta {LB}}_g$: the
category fibred on the category of twisted curves $\wt \MMM_g$
whose fibre on $\stC\to X$ is the
stack $\sta{LB}_{\sta f}$.
\end{rem}

Let $\sta F$ be a line bundle on $\stC$,
whose relative degree is a multiple of $r$; over $X$ we have
$\sta F\to \stC.$
Consider the category ${\sta F^{1/r}}$ of
$r$th roots of $\sta F_S=\sta F\times_X S$
on $\stC_S=\stC\times_X S$ for an $X$-scheme $S$.
More precisely, the category is formed by the
\emph{objects} $(S,\sta M,\sta j)$, where $S$ is
an $X$-scheme, $\sta M$ is a line bundle on
$\stC_S$, and $\sta j$ is an isomorphism
$\sta M^{\otimes r}\xrightarrow{\sim} \sta F_S$.
The \emph{morphisms} $(S,\sta M,\sta j)\to (S',\sta M',\sta j')$ are
pairs
$(m,\sta a)$ as above, with
$\sta a^{\otimes r}$ commuting with $\sta j$ and $\sta j'$.

\begin{rem}
The line bundle $\sta F$ can be regarded as
a section  from $X$ to $\sta{LB}_{\sta f}$.
Consider the $X$-morphism
$\sta k_r\colon \sta{LB}_{\sta f}\to \sta{LB}_{\sta f}$
defined by $(S,\sta M)\mapsto (S,\sta M^{\otimes r})$ and
$(m,\sta a)\mapsto (m,\sta a^{\otimes r})$.
In this way, we can equivalently define ${\sta F^{1/r}}$
as the fibre product
$${\sta F^{1/r}}=(\sta{LB}_{\sta f})\ _{\sta k_r}\!\!\times_{\sta F} X.$$
\end{rem}

\begin{pro}\label{pro:stackLBF/r}
The category ${\sta F^{1/r}}$
satisfies the following properties.
\begin{enumerate}
\item[I.] It is a
Deligne--Mumford stack, \'etale and separated over $X$.
\item[II.] If $\sta F=\Ocal$, it is a group stack on $X$.
\item[III.] For any $\sta F$ the functor
\begin{align*}
\sta m\colon {\Ocal^{1/r}}\times_X {\sta F^{1/r}}&
\to {\sta F^{1/r}}\\
(\sta N,\sta j), (\sta M,\sta k)&
\mapsto (\sta N\otimes \sta M,\sta j\otimes \sta k)
\end{align*}
is an action of ${\Ocal^{1/r}}$ on ${\sta F^{1/r}}$,
and $$(\sta m \times \sta{pr}_2)\colon
{\Ocal^{1/r}}\times_X {\sta F^{1/r}}\to
{\sta F^{1/r}}\times_X {\sta F^{1/r}}$$
is an isomorphism of stacks.
\end{enumerate}\end{pro}
\begin{proof}
Using \cite{Li} or Appendix A,
the proof of the fact that
$\sta F^{1/r}$ is of Deligne--Mumford type only amounts to showing
that the
diagonal is  unramified
(the category of algebraic stacks is closed under fibred products
\cite[4.5]{LM}).
We show that the fibre of the morphism
$\Isom_S(\al,\beta)\to S$ is reduced.
Indeed, at each point $s\colon\Spec k\to S$ the
fibre is either empty, if $\sta M_{\al,s}\ncong \sta M_{\beta,s}$,
or represented by the reduced group scheme $\pmmu_r(k)$, if
$\sta M_{\al,s}\cong \sta M_{\beta,s}$
(the group scheme $\pmmu_r(k)$ acts transitively and freely by
multiplication along the fibres of the line bundle).

The \'etaleness of $\sta F^{1/r}\to X$ claimed in part (I) follows from
the fact that
the relative cotangent complex vanishes.
Indeed, as shown  in \cite[Prop.~3.0.2]{ACV}
and \cite[\S2.1]{AJ},
this is a consequence of the fact that the
relative cotangent complex $\mathbb{L}_{{\sta k}_r}$ vanishes.

In order to show that $\sta F^{1/r}\to X$ is separated
it suffices to consider the case
$\sta F=\Ocal$ and to apply the valuative criterion.
Let $(\Spec R, \sta M_R,\sta j_R)$ be an $r$th root of $\Ocal$
on a twisted curve $\sta C_R$ over a discrete valuation ring $R$.
Over the field of fractions $K$, we assume that $(\Spec K, \sta M_R\otimes K, \sta j_R\otimes K)$ is
trivial: i.e. there exists a trivialization $\sta b_K\colon
\sta M_R\otimes K\to \Ocal$ compatible with $\sta j_R\otimes K$.
Then, the trivialization
$\sta b_K$ can be extended over $R$.
If we assume that $\sta M_R$
is a pullback from $\coa{\sta C_R}$, then we can
focus on
$\coa{\sta C_R}\to \Spec R$ and
the claim follows from the separateness of
the functor of line bundles of degree zero on each irreducible
components of each fibre, \cite{Ra}.
We claim that $\sta M_R$ is a pullback from $\coa{\sta C_R}$, because
the stabilizers at the nodes act trivially on the fibres. This is immediate for
nodes that belong to the closure of the nodes of the generic fibre. Otherwise,
we can focus on a `new' node, whose local picture is
$[\{zy=t\}/\pmmu_l]$ for $t\in (\pi)\subset R$. Note that, here, as soon as
$\pmmu_l$ acts nontrivially on the fibres,
the generic $r$th root of $\Ocal$ is nontrivial (it induces a
nontrivial cyclic covering of the generic fibre).

Part (II) follows from the fact that the stack
${\Ocal^{1/r}}$ is a functor sending
each $X$-scheme $S$ to the groupoid of $r$-torsion line bundles
on $\sta C_S$. In fact such a groupoid is a Picard category:
a symmetric monoid where each object is invertible
and the functors
\begin{align*}
{\Ocal^{1/r}}\times_X{\Ocal^{1/r}}
&\xrightarrow{\  }{\Ocal^{1/r}}&\\
(S,(\sta M_1,\sta j_1),(\sta M_2,\sta j_2))&
\mapsto (S,\sta M_1\otimes \sta M_2,\sta j_1\otimes \sta j_2),\\\\
X&\xrightarrow{\ } {\Ocal^{1/r}} \\
S &\mapsto (S,\Ocal,\id)
\end{align*}and \begin{align*}
{\Ocal^{1/r}}&\xrightarrow{\ } {\Ocal^{1/r}} \\
(S,\sta M,\sta j)&\mapsto (S,\sta M^{\vee},(\sta j^{\vee})^{-1})
\end{align*}
satisfy the law of associativity,
the law of the identity,
and the
law of the inverse.

Part (III) holds because the inverse of the morphism
$\sta m\times \sta{pr}_2$ is given by the functor
$(S,(\sta M_1,\sta k_1),(\sta M_2,\sta k_2))\mapsto
(S,(\sta M_1\otimes \sta M_2^\vee,\sta k_1\otimes(\sta k_2^\vee)^{-1}),
(\sta M_2,\sta k_2))$.
\end{proof}

\begin{defn}[the stack $\ol{\sta F^{1/r}}$]
Note that each object $(S,\sta M, \sta j)$ of $\sta{F}^{1/r}$
over an $X$-scheme $S$ has automorphisms given by multiplication
by an $r$th root of unity along the fibre of $\stL$.
We can eliminate these automorphism by passing to the
corresponding rigidified stack (the general setting is recalled in
\S{}A.2).
The stack $\ol{\sta F^{1/r}}$ is the rigidification of
${\sta F^{1/r}}$ along $\pmmu_r$ in the sense of
\cite{ACV} (see Theorem \ref{thm:rigid}):
$$\ol{\sta F^{1/r}}=({\sta F^{1/r}})^{\pmmu_r}.$$
\end{defn}
\begin{pro}\label{pro:rigidified}
The morphism $\ol{\sta F^{1/r}}\to X$ is
represented by an \'etale and separated $X$-scheme, whose
geometric fibre on
$x$ in $X$ is
$$\{\al\in \Pic (\stC_x)\mid \al^{r}=[\sta F_x]\}.$$
Assume that $\sta F$ is the trivial line bundle $\Ocal$ of $\stC$; then,
$\ol{\Ocal^{1/r}}$ is a group scheme.
\end{pro}
\begin{proof}
The morphism ${\sta F^{1/r}}\to X$
is \'etale and separated, Proposition
\ref{pro:stackLBF/r}.
The morphism
${\sta F^{1/r}}\to \ol{\sta F^{1/r}},$
is surjective and indeed
a morphism
``locally isomorphic''
to $\sta B\pmmu_r$ on $\ol{\sta F^{1/r}}$
(see  \cite[I.~Prop.~3.0.2, (2)]{Ro} and \cite[Thm.~5.1.5]{ACV}).
Therefore, $\ol{\sta F^{1/r}}\to X$ is \'etale and separated.

Finally $\ol{\sta F^{1/r}} \to X$ is representable
because for any algebraically closed field $k$
and any $\tau\in \ol{\sta F^{1/r}}(\Spec k)$
the automorphism group of $\tau$ is trivial (this condition suffices by
a criterion due to Abramovich and Vistoli,
\cite[Lem.~4.4.3]{AV}).
Indeed,
the objects of
$\ol{\sta F^{1/r}}=({\sta F^{1/r}})^{\pmmu_r}$ over
an algebraically closed field $k$
are the same as the objects of ${\sta F^{1/r}}$ over $k$.
Now, the automorphism group of the
$\pmmu_r$-rigidified object $\tau$
is the quotient by $\pmmu_r(k)$ of the automorphism
group of $\tau$ in ${\sta F^{1/r}}$, Theorem \ref{thm:rigid}, (3).
This quotient can be easily seen to be trivial: note that $\pmmu_r(k)\to \Aut(\tau)$ is surjective,
because  an automorphism
of an $r$th root
$(\sta M,\sta j\colon \sta{M}^{\otimes r}\xrightarrow{\sim} \sta F)$
on a twisted curve $\stC\to \Spec k$
$$\sta a\colon \sta M\xrightarrow{\sim} \sta M \quad \quad
\sta a^{\otimes r}\circ \sta j=\sta j$$
necessarily satisfies $\sta a^{\otimes r}=\id$ and, since $\sta C\otimes k$ is connected,
$\sta a$ is given by the multiplication by an $r$th root
of unity along
the fibre of $\sta M$.

The rest of the proposition follows immediately from
Proposition \ref{pro:stackLBF/r} (the group structure
and the action of ${\Ocal^{1/r}}$
on ${\sta F^{1/r}}$
descend to $\ol{\sta F^{1/r}}$ and
their relations are preserved).
\end{proof}
\begin{pro} \label{pro:finiteness}
Assume that $\sta C\to X$ is a twisted curve and $\sta F\to \sta C$
is a line bundle  satisfying the following conditions:
for any geometric point
$x$ in $X$ the line bundle $\sta F_x$ has $r^{2g}$ $r$th roots on
$\stC_x$ up to isomorphism.
Then, we have the following properties.
\begin{enumerate}
\item[(1)] $\ol{\Ocal^{1/r}}$ is a finite group scheme
and $\ol{\sta F^{1/r}}$
is a finite torsor under $\ol{\Ocal^{1/r}}$.
\item[(1')]
${\Ocal^{1/r}}$ is a finite group stack
and ${\sta F^{1/r}}$ is a finite torsor under
${\Ocal^{1/r}}$.
\end{enumerate}
\end{pro}
\begin{proof}
The morphism $p\colon Y=\ol{\sta F^{1/r}}\to X$
is \'etale and the geometric fibres
are reduced and consist
of $n=r^{2g}$ distinct points. This condition is
sufficient for the properness of $p$: indeed it implies that
$p$ is an isomorphism when $n$ equals $1$,
and it implies properness by induction on the degree
(we can check properness after an \'etale base change,
and we note that $Y\times_X Y$ is finite on $Y$
because $(Y\times_X Y)$ minus the diagonal $\Delta(X)$
is \'etale and separated on $Y$,
and has constant reduced fibre
consisting of $n-1$ points).

Now, note that if $\sta F_x$ has $r^{2g}$ distinct $r$th roots on
$\sta C_x$, then
the structure sheaf of $\stC_x$ also has $r^{2g}$ distinct roots.
Therefore, $\ol{\Ocal^{1/r}}$ is a
finite group and $\ol{\sta F^{1/r}}$ is a finite torsor
under $\ol{\Ocal^{1/r}}$, because
the geometric fibres of $\ol{\sta F^{1/r}}\to X $ are nonempty,
and $\ol{\Ocal^{1/r}}\times _X \ol{\sta F^{1/r}}\to \ol{\sta F^{1/r}}\times_X\ol{\sta F^{1/r}}$
is an isomorphism.
Finally the point $(1')$ follows from Proposition \ref{pro:stackLBF/r},
and the fact that $\sta F^{1/r}$
is proper on $\ol{\sta F^{1/r}}$.
\end{proof}

\subsection{The Kummer theory of a twisted curve}
\label{sect:kummer}
Since, by the above proposition, ${\sta F^{1/r}}$ is finite as soon as the fibre is constant,
we focus on the geometric fibres of $\ol{\sta F^{1/r}}\to X$.
\begin{notn}[multiplicity of a line bundle at a smooth stack-theoretic point]
We consider $1$-dimensional
Deligne--Mumford stacks $\sta X$ with nodal singularities
and trivial stabilizers except for a finite number of points.
As noted in Remark \ref{rem:cangener}
for each smooth point $\sta p$
there is a
canonical embedding $$\sta j\colon \sta B(\Aut(\sta p))\to \stC$$
and
a canonical generator of
$\Pic \sta B(\Aut(\sta p))$ induced by the
$\Aut(\sta p)$-linearized tangent space $\sta T$.
With such a canonical generator,
the group $\Pic \sta B(\Aut(\sta p))$
can be regarded as $\ZZ/\size {\Aut(\sta p)}\ZZ$.
Then for any line bundle
$\sta F$ at any smooth point $\sta p$ of $\sta X$
we define a
multiplicity index
$$0\le \mult _\sta p \sta F\le \size {\Aut(\sta p)}-1$$
via the pullback homomorphism
$$\Pic \sta X \xrightarrow{\sta j^*}\Pic \sta B(\Aut(\sta p)).$$
If $\sta e$ is a node the notation $\mult_{\sta e}\sta F$ does not make sense,
because there is no canonical generator for $\Pic \sta B(\Aut(\sta p))$.
Instead, we will consider the pullback of $\sta F$ on the normalization,
which is  smooth, $1$-dimensional, and contains two distinct geometric points
$\sta p_1$ and $\sta p_2$ lying over the node $\sta e$. By a slight abuse of notation
we shall write $\mult_{\sta p_1}\sta F$ and $\mult_{\sta p_2}\sta F$
for the multiplicity index at $\sta p_1$ and $\sta p_2$ of the pullback of $\sta F$ on the normalization.
\end{notn}
\begin{theor}\label{thm:rootsnum}
Let $\stC\to \Spec k$
be a twisted curve with stabilizers of order $l(\sta e)$
on each node $\sta e$.
Let $\sta F$ be a line bundle on
$\stC$  whose total degree
is a multiple of $r$.

The number  of
$r$th roots of $\sta F$ is $r^{2g}$ if and only if
the two following numerical conditions are satisfied.
First, for each nonseparating node $\sta e$, $r$ divides
$$l(\sta e), \; m_1(\sta e), \text{ and } m_2(\sta e)$$
where the normalization
$\stC$ at $\sta e$ is the connected stack $\stC_0$,
$\sta p_1$, $\sta p_2\in \stC_0$ map to $\sta e$, and
$m_i(\sta e)$ equals $\mult_{\sta p_i}(\sta F)$.
Second, for each separating node $\sta e$, $r$ divides
$$d_1(\sta e)l({\sta e}) \text{ and } d_2(\sta e)l(\sta e)$$
where the normalization of $\stC$ at $\sta e$
is $\stC_1\sqcup \stC_2$
and $d_i\in ({1}/{l(\sta e)})\ZZ$ is the total degree of $\sta F$
on the curve $\stC_i$.

In particular, for any line bundle
on $\coa{\stC}$, whose total degree is a multiple of $r$,
there are exactly $r^{2g}$ roots of the pullback on $\stC$ if and only if
$r$ divides $l(\sta e)$ for each nonseparating node, and
$r$ divides $l(\sta e) d_i(\sta e)$
for each separating node.
\end{theor}

We prove the above criterion at the end of the section (\ref{point}).
The main tool is the exactness of the
following diagram (Figure 3.1) and its interpretation
for twisted curves (\ref{point}, Figure 3.2).
\begin{theor}\label{thm:roots}
Let $\pi \colon \sta X\to \coa{\sta X}$ be a
tame and proper Deligne--Mumford stack
of dimension $1$ with nodal singularities and assume that
$\Aut(\sta p)$ is trivial
except for a finite number of points
$\sta p_i\colon \Spec k\to \sta X$. We
write  $\Gamma_i=\Aut(\sta p_i)$
and we denote by
$\sta j$ the canonical embedding
$\sqcup_i \sta B\Gamma_i\to \sta X$.

Then, the direct images $R\pi_*\pmmu_r$ and
$R\pi_*\GG_m$ are represented by
complexes on $\coa{\sta X}$
fitting in the exact sequences
\begin{align}
1\to \pmmu_r\to R\pi_*\pmmu_r&\to R\pi_*\pmmu_r/\pmmu_r\to 1,
\label{eq:hypermu}\\
1\to \GG_m\to R\pi_*\GG_m&\to R\pi_*\GG_m/\GG_m\to 1 \label{eq:hypergm}.
\end{align}
The long exact sequence of hypercohomology
of \eqref{eq:hypermu}
yields
\begin{multline*}
1 \to H^1(\coa{\sta X},\pmmu_r)\to H^1({\sta X},\pmmu_r)\to \textstyle{\prod_i H^1(\Gamma_1,\pmmu_r)}
  \\\xrightarrow{\ \ \delta\ \ } H^2(\coa{\sta X},\pmmu_r)\to H^2({\sta X},\pmmu_r)\to
  \textstyle{\prod_i H^2(\Gamma_1,\pmmu_r)}\to 1,
\end{multline*}
and the long exact sequence of hypercohomology
of \eqref{eq:hypergm} in degree $1$ yields
\begin{align*}
1\to \Pic \coa{\sta X}\xrightarrow {\pi^*} \Pic{\sta X}\to
\textstyle{\prod_i \Pic(\sta B\Gamma_i)}\to 1
\end{align*}
Both exact sequences above fit in the
following commutative diagram (see the
second and the third columns), where
all vertical
and horizontal sequences are exact

\begin{figure}[h]\label{fig:thm}
\begin{equation*}
{ \xymatrix@C=0.25cm{
\Pic {\coa{\sta X}}\ar[r]&
H^2({\coa{\sta X}},\pmmu_r)\ar[r]&
1\ar[r] &
1\ar[r] & 1\ar[r] &1\\
1\ar[r]\ar[u] &
{\textstyle\prod H^1(\Gamma_i,\pmmu_r)}\ar[r]\ar[u]_{\delta}&
{\textstyle\prod\Pic(\sta B\Gamma_i)}\ar[r]^r\ar[u]&
{\textstyle \prod\Pic(\sta B\Gamma_i)}
                   \ar[r]\ar[u]&
{\textstyle \prod H^2(\Gamma_i,\pmmu_r)}\ar[u]\ar[r]&1\ar[u]\\
1\ar[r]\ar[u] &H^1(\sta X,\pmmu_r)     \ar[r]\ar[u]&\Pic\sta X
\ar[r]^r\ar[u]_{\sta j^*} &\Pic\sta X \ar[r]
\ar[u]_{\sta j^*}&H^2(\sta X, \pmmu_r)\ar[u]\ar[r]&1 \ar[u]\\
1\ar[r]\ar[u] &H^1({\coa{\sta X}},\pmmu_r)     \ar[r]\ar[u] &\Pic {\coa{\sta X}} \ar[r]^r\ar[u]_{\pi^*}
&\Pic {\coa{\sta X}} \ar[r]\ar[u]_{\pi^*} &H^2({\coa{\sta X}},\pmmu_r)\ar[u]\ar[r]&1 \ar[u]\\
1\ar[u]\ar[r]& 1\ar[u]\ar[r]&1\ar[u]\ar[r]&
1\ar[u]\ar[r] &
{\textstyle\prod H^1(\Gamma_i,\pmmu_r)}\ar[u]_{\delta}\ar[r]&{\textstyle\prod\Pic(\sta B\Gamma_i)}.\ar[u]
}}
\end{equation*}
\begin{caption}{the diagram can be regarded as
a double complex $K^{\bullet,\bullet}$
periodic of period $(3,-3)$.}
\end{caption}
\end{figure}
\end{theor}


\begin{proof}
The exactness of
$
1\to \pmmu_r\to \GG_m
\to \GG_m\to 1$
holds in the \'etale topology of $\sta X$
by the same argument of
\cite[II.2.18b]{Mi}. Our goal is to compare the long exact sequences of
cohomology of the Kummer sequence of $\sta X$ and of $\coa{\sta X}$.

We show that there are three complexes on
${\coa{\sta X}}$ representing the direct images
in the derived category $R\pi_* \pmmu_r$, $R\pi_* \GG_m$, and $R\pi_* \GG_m$
and fitting in a short exact sequence as follows
\begin{equation}\label{eq:Rpi}
1\to R\pi_* \pmmu_r\to
R\pi_*\GG_m\to R\pi_*\GG_m\to 1.
\end{equation}
Indeed, in the category of sheaves of abelian groups
on $\sta {\coa{\sta X}}$, we can consider
three injective resolutions
\begin{align*}
\pmmu_r\to I_1^\bullet,\\
\GG_m \to I_2^\bullet,\\
\GG_m \to I_3^\bullet,
\end{align*}
fitting in
\[ \xymatrix@R=0.4cm{
1\ar[r] &I_1^\bullet \ar[r] &I_2^\bullet                 \ar[r] &I_3^\bullet\ar[r]& 1\\
1\ar[r] &\pmmu_r     \ar[r]\ar[u] &\mathcal \GG_m
\ar[r]^r\ar[u] &\mathcal \GG_m \ar[r]\ar[u]& 1,
} \]
where all the
horizontal sequences are exact (note that this happens because the
category of sheaves of abelian groups has enough injectives, and, given the injective resolutions
$I_1^\bullet$ and $I_3^\bullet$, we can also
find a third injective resolution $I_2^\bullet$ making all horizontal sequences
exact). Now, applying
the direct image via
$\pi\colon \sta X\to {\coa{\sta X}}$ to $1\to I_1^\bullet\to I_2^\bullet\to I_3^\bullet\to 1$,
we get
the short exact sequence
\eqref{eq:Rpi}, because all the sheaves involved  are injective.

The sequence \eqref{eq:Rpi} allows the comparison of the long
exact cohomology sequences attached to the
Kummer sequences of $\sta X$ and of $\coa {\sta X}$. Indeed,
we show that the Kummer sequence of $\coa{\sta X}$
injects in \eqref{eq:Rpi}.
We have
$\pi_*\pmmu_r=\pmmu_r$ and
$\pi_*\GG_m=\GG_m$, because
$\pi_*\Ocal_{\stC}$ is canonically isomorphic to $\Ocal_{\coa{\stC}}$
by Keel and Mori's theorem \cite{KM}, and the isomorphism yields the identity
$\pi_*\pmmu_r=\pmmu_r$ because it identifies the $r$th roots of unity.
We deduce that, on ${\coa{\sta X}}$,  the sequence $1\to \pmmu_r\to \GG_m
\to \GG_m\to 1$ injects into \eqref{eq:Rpi}. Therefore, we write
\begin{equation}\label{eq:diagr}
\xymatrix@R=0.4cm{
1\ar[r] &R\pi_*\pmmu_r     \ar[r] &R\pi_*\GG_m\ar[r]^r
&R\pi_*\GG_m\ar[r]& 1\\
1\ar[r] &\pmmu_r     \ar[r]\ar[u] &\GG_m \ar[r]^r\ar[u] &\GG_m\ar[r]\ar[u]& 1\\
& 1\ar[u]&1\ar[u]&1\ar[u]&.
}
\end{equation}
Note that, each vertical sequence
in the diagram above is of the form
$\pi_*\sta J\to R\pi_*\sta J$, where, on the one hand,
$\pi_*\sta J$ is concentrated in degree $0$:
$(1\to A^0\to 1\to 1\to \dots)$,
on the other hand $R\pi_*\sta J$ is
$(1\to B^0\to B^1\to B^2\dots )$
and the morphism of complexes is given by
\begin{equation*}
\label{eq:diag}
\xymatrix@R=0.4cm{
1\ar[r] &B^0\ar[r] &B^1\ar[r] &B^2 \\
1\ar[r] &A^0\ar[r] \ar[u] &1\ar[u]\ar[r] &1\ar[u] &, \\
}
\end{equation*}
where $A^0\to B^0$ is injective. By taking quotients, we get
complexes representing $R\pi_*\sta J/\pi_*\sta J$ for $\sta
J=\pmmu_r$ or $\GG_m$. By the exactness of the horizontal
sequences in \eqref{eq:diagr}, we get the following
diagram where horizontal and vertical sequences are exact.
\begin{equation}\label{eq:3by3}
\xymatrix@R=0.4cm{
&1&1& 1\\
1\ar[r] &R\pi_*\pmmu_r/\pmmu_r\ar[u] \ar[r] &R\pi_*\GG_m/\GG_m\ar[u]\ar[r]
&R\pi_*\GG_m/\GG_m\ar[u]\ar[r]& 1\\
1\ar[r] &R\pi_*\pmmu_r \ar[u]    \ar[r] &R\pi_*\GG_m\ar[u]\ar[r]^r
&R\pi_*\GG_m\ar[u]\ar[r]& 1\\
1\ar[r] &\pmmu_r     \ar[r]\ar[u] &\GG_m \ar[r]^r\ar[u] &\GG_m\ar[r]\ar[u]& 1\\
& 1\ar[u]&1\ar[u]&1\ar[u]&
}
\end{equation}

We calculate the hypercohomology groups of ${\coa{\sta X}}$ with
respect to the complexes $\GG_m,$ $\pmmu_r,$ $R\pi_*\GG_m$, $R\pi_*\pmmu_r$,
$R\pi_*\GG_m/\GG_m$, and $R\pi_*\pmmu_r/\pmmu_r$
in the diagram above.
Since the hypercohomology groups coincide with
the cohomology groups when the complex is
concentrated in degree $0$, for any sheaf $E$ on ${\coa{\sta X}}$ we have
\begin{equation}
\label{eq:hyperX}
\HH^j({\coa{\sta X}}, E)=H^j({\coa{\sta X}},E).
\end{equation}
By the Leray spectral sequence, for any sheaf of groups $\sta J$
on $\sta X$, we have
\begin{equation}\label{eq:hyperstX}
\HH^j({\coa{\sta X}},R\pi_*\sta J)=H^j(\sta X,\sta J).
\end{equation}
Finally the spectral sequence
$E^{p,q}_2=H^p({\coa{\sta X}},H^q(R\pi_*\sta J/ \pi_*\sta J))$ abuts to
 $\HH^{p+q}({\coa{\sta X}}, R\pi_*\sta J/ \pi_*\sta J)$.
The sheaf  $H^q(R\pi_*\sta J/ \pi_*\sta J)$
is equal to $R^q\pi_*\sta J$ for $q>0$, and vanishes
otherwise.
We now notice that $R^q\pi_*\sta J$ is supported on the
points with nontrivial stabilizer. Indeed, it is shown in
\cite[Prop.~A.0.1]{ACV} that the stalk
 of $R^q\pi_*\sta J$ at a point
 $\coa{\sta p}\in {\coa{\sta X}}$
is canonically isomorphic to the $q$th cohomology group
$H^q(\Aut(\sta p),\sta J_{\sta p})$.
Therefore, we have
\begin{equation}\label{eq:hypersp}
\HH^j({\coa{\sta X}}, R\pi_*\sta J/ \pi_*\sta J)=\begin{cases}
0 & \text{if $j=0$}\\
\textstyle{\prod_{i}H^j(\Gamma_i,\sta J_{\sta p_i})}& \text{if  $j>0$,}
\end{cases}
\end{equation}
because $R^q\pi_*\sta J$ is supported on
the points $\coa{\sta p_i}$, and the differentials of
$E^{p,q}_2$ vanish.

Passing to the long exact hypercohomology sequences in all directions
in the diagram \eqref{eq:3by3}, we get
the following diagram. The fact that
the diagram \eqref{eq:3by3} commutes
implies that all the squares
below commute except for the squares involving the boundary homomorphisms,
which anticommute. However, note that
all the composition of homomorphism in the
squares involving the boundary homomorphisms
are trivial, so the following diagram
is indeed commutative
\begin{equation*}
{ \xymatrix@C=0.30cm{
H^1({\coa{\sta X}},\GG_m)\ar[r]&
H^2({\coa{\sta X}},\pmmu_r)\ar[r]&
1\ar[r]&
1\ar[r]& 1\ar[r]& 1
\\
1\ar[r]\ar[u] &
{\textstyle\prod H^1(\Gamma_i,\pmmu_r)}\ar[r]\ar[u]^{\delta}&
{\textstyle\prod H^1(\Gamma_i,\GG_m)}\ar[r]\ar[u]&
{\textstyle\prod H^1(\Gamma_i,\GG_m)}
                   \ar[r]\ar[u]&
{\textstyle \prod H^2(\Gamma_i,\pmmu_r)}\ar[u]\ar[r]&1\ar[u]
\\
1\ar[r]\ar[u] &H^1(\sta X,\pmmu_r)     \ar[r]\ar[u]&H^1(\sta X,\GG_m)
\ar[r]\ar[u] &H^1(\sta X,\GG_m)\ar[r]
\ar[u]&H^2(\sta X, \pmmu_r)\ar[u]\ar[r]&1\ar[u]
\\
1\ar[r]\ar[u] &H^1({\coa{\sta X}},\pmmu_r)     \ar[r]\ar[u] &H^1({\coa{\sta X}},\GG_m)\ar[r]\ar[u]
&H^1({\coa{\sta X}},\GG_m) \ar[r]\ar[u]
&H^2({\coa{\sta X}},\pmmu_r)\ar[u]\ar[r]&1\ar[u]
\\
1\ar[u]\ar[r] & 1\ar[u]\ar[r]&1\ar[u]\ar[r]&1\ar[u]\ar[r]&
{\textstyle\prod H^1(\Gamma_i,\pmmu_r)}\ar[u]^{\delta}\ar[r]&{\textstyle\prod H^1(\Gamma_i,\GG_m)}.\ar[u]
}}
\end{equation*}
In the diagram above, we used the following facts.
\begin{enumerate}
\item
For any proper stack $\sta Y$
the $r$th power homomorphism is surjective in $H^0(\sta Y,\GG_m)$.
\item
The coarse space ${\coa{\sta X}}$ is a $1$-dimensional scheme with nodal singularities, therefore
we have $H^2({\coa{\sta X}},\GG_m)=1$.
\item
For $j>0$, the cohomology groups of $\Gamma_i=\Aut(\sta p_i)$ with coefficient in
$(\GG_m)_{\sta p_i}=\GG_m$ with a trivial $\Gamma_i$-action
are given by
$$H^j(\Gamma_i,\GG_m)=\begin{cases}1 & \text {if $j$ is even,} \\
\ZZ/l(i)\ZZ &\text{if $j$ is odd,}\end{cases}$$
where $l(i)=\size{\Gamma_i}$.
\end{enumerate}

After the identifications
of $H^1(\sta Y,\GG_m)$ with $\Pic(\sta Y)$,
this shows the exactness of the diagram of
Figure 3.1.
\end{proof}
\begin{rem}[the boundary homomorphism
$\delta$]
We show that the coboundary homomorphism
$$\delta\colon \textstyle{\prod_i H^1(\Gamma_i,\pmmu_r)} \to H^2({\coa{\sta X}},\pmmu_r)$$
which maps vertically in the commutative diagram
of Figure 3.1 can be explicitly computed.
When $\sta X$ is a twisted curve,
this allows to provide the explicit description of
the exact sequence
\begin{multline}\label{eq:longexa}
1\to H^1({\coa{\sta X}},\pmmu_r)\to H^1(\sta X,\pmmu_r)\to
\textstyle{\prod H^1(\Gamma_i},\pmmu_r)\\\xrightarrow{\delta}
H^2({\coa{\sta X}},\pmmu_r)\to H^2(\sta X,\pmmu_2)\to
\textstyle{\prod H^1(\Gamma_i},\pmmu_r)\to 1
\end{multline}
and to prove Corollary \ref{cor:root_tw} and the weaker version
Corollary \ref{cor:rdividesl} stated above.

First, for any $1$-dimensional stack $\sta X$
satisfying the hypotheses of the theorem, we
use the double complex of Figure 3.1 to
define a $\ZZ\times \ZZ$-graded
complex $K^{\bullet,\bullet}$ as follows.
\begin{equation*}
{ \xymatrix@C=0.45cm{
&1&1& 1& 1& \\
1\ar[r] &
{\textstyle\prod_iH^1(\Gamma_i,\pmmu_r)}\ar[r]^{d'}\ar[u]&
{\textstyle\prod_iH^1(\Gamma_i,\GG_m)}\ar[r]^{d'}\ar[u]&
{\textstyle\prod_iH^1(\Gamma_i,\GG_m)}
                   \ar[r]^{d'}\ar[u]&
{\textstyle \prod_{i}H^2(\Gamma_i,\pmmu_r)}\ar[u]\ar[r]&1\\
1\ar[r] &H^1(\sta X,\pmmu_r)     \ar[r]^{d'}\ar[u]^{d''}&H^1(\sta X,\GG_m)
\ar[r]^{d'}\ar[u]^{d''} &H^1(\sta X,\GG_m)\ar[r]^{d'}
\ar[u]^{d''}&H^2(\sta X, \pmmu_r)\ar[u]^{d''}\ar[r]&1 \\
1\ar[r] &H^1({\coa{\sta X}},\pmmu_r)
 \ar[r]^{d'}\ar[u]^{d''} &H^1({\coa{\sta X}},\GG_m)\ar[r]^{d'}\ar[u]^{d''}
&H^1({\coa{\sta X}},\GG_m) \ar[r]^{d'}\ar[u]^{d''}
&H^2({\coa{\sta X}},\pmmu_r)\ar[u]^{d''}\ar[r]&1 \\
 & 1\ar[u]&1\ar[u]&1\ar[u]&1\ar[u]&,
}}
\end{equation*}
Here $K^{0,0}$ equals $H^1({\coa{\sta X}},\pmmu_r)$
and all terms $K^{p,q}$ coincide with
the terms of
the diagram of Figure 3.1 appearing in Theorem \ref{thm:roots}
except for $K^{3,-1}=1$ and $K^{0,3}=1$.
We also
set the differentials $d'(p,q)\colon K^{p,q}\to K^{p+1,q}$ and
$d''(p,q) \colon K^{p,q}\to K^{p,q+1}$
as in Figure 3.1, except for
$d''(0,2)$ and $d''(3,-1)$ which are trivial.

\vspace{0.2cm}
\noindent Now, by means of the double complex $K^{\bullet, \bullet}$ we
identify $\delta$ with a differential $'\!d_3$ of the spectral sequence
associated to the filtration
$'K_p=\sum_{i\ge p}K^{i,j}$.
This amounts to the following \emph{explicit description} of $\delta$.
For any $x\in \prod_iH^1(\Gamma_i,\pmmu_r)$,
$\delta(x)$ is the (unique) element of $H^2({\coa{\sta X}}, \pmmu_r)$
for which there exist $y_1\in H^1(\sta X, \GG_m)$ and
$y_2\in H^1({\coa{\sta X}}, \GG_m)$ satisfying
\begin{equation}\label{eq:syst}
\begin{cases}
d'(x)=d''(y_1), \quad \quad\\
d'(y_1)=d''(y_2), \quad \quad\\
d'(y_2)= \delta(x).
\end{cases}\end{equation}

\vspace{0.2cm}

Indeed, we consider the spectral sequences $'\!E_\bullet$ and $''\!E_\bullet$
with
$$'\!E_2^{p,q}=H_{d'}^p(H_{d''}^q(K^{\bullet,\bullet})) \quad \quad \quad
''\!E_2^{p,q}=H_{d''}^q(H_{d'}^p(K^{\bullet,\bullet}))$$
associated to the filtrations
$'K_p=\sum_{i\ge p}K^{i,j}$
and $''K_q=\sum_{j\ge q}K^{i,j}$.
Note that $d'$ is exact everywhere; therefore,
$H_{d'}(K^{\bullet,\bullet})$
vanishes identically.
So, $''\!E_2^{p,q}$ is zero; hence,
also $'\!E_h^{p,q}$ vanishes for all $p$ and $q$ and
for $h$ sufficiently large. Since $'\!E_h^{p,q}$
is constant for $h\ge 4$, the homomorphism
\begin{equation}\label{eq:d3}
\,'\!E^{0,2}_3=\coker (d''(0,1))
\xrightarrow{\ \ \ '\!d_3\ \ \ }
\,'\!E^{3,0}_3=\ker (d''(3,0)).\quad \quad \end{equation}
is an isomorphism.
In fact, by construction,
the coboundary homomorphism $\delta$ coincides with
the composite of
$K^{0,2}\to \coker (d''(0,1))$
and $'\!d_3$. The definition of the differential
$'\!d_3$ for the spectral sequence $'\!E_3^{p,q}$  yields
the equations \eqref{eq:syst} above.
\end{rem}

In the following corollary we apply
Theorem \ref{thm:roots}
and the explicit computation of
the boundary homomorphism $\delta$ given above
to the case of a twisted curve.
See also \ref{point}, Figure 3.2 where we summarize the result.
\begin{cor}\label{cor:root_tw}
Let $\pi\colon\stC\to {\coa{\sta C}}$ be a twisted curve. Write
$\Gamma_{\sta e}$ for $\Aut(\sta e)$ for any node $\sta e\in E$
and $l(\sta e)$ for its order.
We have a canonical isomorphisms
\begin{equation}\label{eq:h2}
H^2(\stC,\pmmu_r)\cong H^2({\coa{\sta C}},\pmmu_r)\cong(\ZZ/r\ZZ)^V.
\end{equation}
Choose an orientation of the dual graph of $\stC$; then,
we have a canonical isomorphism identifying
$H^1(\Gamma_{\sta e},\pmmu_r)$ with the $r$-torsion of
$\ZZ/l(\sta e)\ZZ$:
\begin{equation}\label{eq:h1}
H^1(\Gamma_{\sta e},\pmmu_r)\cong
(l(\sta e)/\hcf\{l(\sta e), r\})\ZZ/l(\sta e)\ZZ \quad  \forall \sta e.
\end{equation}
With respect to these identifications,
$\delta\colon
\prod_{\sta e} H^1(\Gamma_{\sta e},\pmmu_r)
\to H^2({\coa{\sta C}},\pmmu_r)$
in Figure 3.1
as the composite homomorphism
of the boundary of the
chain complex $\Ccal_\bullet(\Lambda, \ZZ/r\ZZ)$
of the dual graph $\Lambda$ of $\stC$
and of the homomorphism
$\prod_{\sta e} H^1(\Gamma_{\sta e},\pmmu_r)
\to (\ZZ/r\ZZ)^E$
given by the product over $E$ of
$$(l(\sta e)/\hcf\{l(\sta e), r\})\ZZ/l(\sta e)\ZZ
\xrightarrow{\ \times (r/l({\sta e}))\ }
(r/ \hcf\{l(\sta e), r\})\ZZ/r\ZZ \into
\ZZ/r\ZZ
\quad  \forall \sta e.$$
\end{cor}
\begin{proof}
First, we see the preliminary case
of the stack $C[p/l]$. Then, we normalize the twisted curve  and use the
presentation of smooth stacks of dimension 1 given in
Section \ref{exa:onedim}.
\begin{lem}\label{lem:onesmooth}
Let  $C$ be a proper, connected, smooth, and reduced curve.
For $p\in C$ and $l$ an invertible integer, consider the stack $C[p/l]$:
the point $C[p/l]$ lying
over $p\in C$ has automorphism group $\Gamma$ of order $l$.
The boundary homomorphism $\delta$ of
Figure 3.1
is the composite homomorphism
\begin{multline*}H^1(\sta B \Gamma,\pmmu_r)\cong(l/\hcf\{l, r\})\ZZ/l\ZZ
\xrightarrow{\ \times (r/l)\ }
(r/\hcf\{l, r\})\ZZ/r\ZZ \into
\ZZ/r\ZZ\cong H^2(C,\pmmu_r),\end{multline*}
where we used the canonical isomorphisms
$\Pic (\sta B\Gamma)\cong \ZZ/l\ZZ$
(Remark \ref{rem:cangener})
and $H^2(C,\pmmu_r)\cong\ZZ/r\ZZ$.
\end{lem}
\begin{proof}
As in Remark \ref{rem:cangener}, we denote by $\sta T$ the
tangent space at $\sta p$.
Now, we apply to $\sta T^{\otimes m}\in H^1(\Gamma,\pmmu_r)$
the definition of $\delta$ given in
\eqref{eq:syst}. Let $\sta M$ be a line bundle
on $C[p/l]$ such that $\sta j^*\sta M=\sta T^{\otimes m}$; let
$A$ be a line bundle on $C$, such that $\pi^*A=\sta M^{\otimes r}$.
Then, $\delta(\sta T^{\otimes m})$ is $\deg A \mod r$.
Finally, by Proposition \ref{pro:deg}, we have
\begin{align*}
\deg(A)=\deg(\sta M^{\otimes r})=r\deg(\sta M)
       =r(k+m/l)\label{eq:why}=rk+mr/l.
\end{align*}
This is indeed the claim of the lemma.
\end{proof}

Now let $\sta D\to \stC$ and
${\coa{\sta D}}\to {\coa{\sta C}}$ be the normalization
of $\stC$ and ${\coa{\sta C}}$.
Note that $\sta D$ has trivial stabilizers
except for a finite set of smooth points $I$.
Therefore, Theorem \ref{thm:roots} applies.
Furthermore, for any point $\sta p\in \sta D$
with $\textstyle{|}\Aut(\sta p)\textstyle{|}=l>1$,
there is a natural projection
\begin{equation}\label{eq:Dproj}
\sta D\to {\coa{\sta D}}\big[\coa{\sta p}/l\big].
\end{equation}
Using Lemma \ref{lem:onesmooth},
we get a canonical isomorphism from
$H^1(\Aut(\sta p),\pmmu_r)$ to $(l/\hcf\{l, r\})\ZZ/l\ZZ$
and a homomorphism
\begin{equation}\label{eq:onepoint}
H^1(\Aut(\sta p),\pmmu_r)\xrightarrow{\ \times (r/l)\ }
H^2({\coa{\sta D}},\pmmu_r).
\end{equation}
Note that the morphism
$\delta_{\sta D}\colon
\prod_{\sta p\in I}H^1(\Aut(\sta p),\pmmu_r)\to H^2({\coa{\sta D}},\pmmu_r)$
coincides with \eqref{eq:onepoint} on each factor,
because $\delta$ commutes with  pullbacks
via the projections \eqref{eq:Dproj}.

Now, consider $\sta D\to \stC$.
For each node $\sta e$ of $\stC$
there are two smooth points
$\sta p_{\sta e, +}, \sta p_{\sta e, -}\in \sta D$
over $\sta e$ with stabilizer $\Gamma_{\sta e}$,
we denote them according to the orientation of the dual graph.
Now, $H^1(\Gamma_{\sta e},\pmmu_r)$
is the $r$-torsion subgroup of
$H^1(\Gamma_{\sta e},\GG_m)$, which is canonically generated by
the tangent space $\sta T_{\sta e,+}$ at
the point of $\sta p_{\sta e,+}\in \sta D$.
Thus, by raising $\sta T_{\sta e,+}$ to the
$(l(\sta e)/\hcf\{l(\sta e), r\})$th power, we get the canonical
isomorphism \eqref{eq:h1}.

The cohomology homomorphisms
induced by pullback via $\sta D\to \sta C$
are compatible with the diagrams
of Figure 3.1 for $\sta C$ and $\sta D$.
Now,
$$\delta_{\stC}\colon \textstyle{\prod_{\sta e}}
H^1(\Gamma_{\sta e},\pmmu_r)
\to H^2({\coa{\sta C}},\pmmu_r)$$
fits in the following commutative diagram
\[ \xymatrix@R=1cm{
 \textstyle{\prod_{\sta e}}
H^1(\Gamma_{\sta e},\pmmu_r)
\ar[r]\ar[r]^{\delta_{\stC}}\ar[d]_{\prod _{\sta e}g_{\sta e}}&
H^2({\coa{\sta C}},\pmmu_r)\ar[d]^\cong \\
\prod_{\sta e, i=+,-}
H^1(\Aut(\sta p_{\sta e,i}),\pmmu_r)
\ar[r]^{\ \ \ \ \ \ \ \ \ \delta_{\sta D}}&  H^2({\coa{\sta D}},\pmmu_r),
} \]
where the homomorphism induced via pullback
$H^2(\coa{\sta D},\pmmu_r)\to H^2(\coa{\sta C},\pmmu_r)$
is invertible and
the morphisms
$g_{\sta e}$ are induced
by pullback via $\sta D\to \sta C$
and map as
$$\quad \quad \sta T_{\sta e,+}\to (\sta T_{\sta e,+},
(\sta T_{\sta e,+})^\vee) \quad \quad \forall \sta e.$$
This happens
because, by the definition of twisted curves \ref{defn:twisted}, the local
picture of
$\stC$ at $\sta e$ is $\{zw=0\}$ with
$\pmmu_{l(\sta e)}$ acting as $(z,w)\mapsto (\xi_{l(\sta e)}z,\xi_{l(\sta e)}^{-1} w)$.
\end{proof}

\begin{point}\label{point}
Proof of Theorem \ref{thm:rootsnum}.\emph{
We prove the claim by chasing in the following diagram derived from
Theorem \ref{thm:roots} and Corollary \ref{cor:root_tw}. }\end{point}
{\
\begin{figure}[!h]\label{fig:bigdiag}
\begin{equation*}
{\xymatrix@C=0.5cm{
P\ar[r]&
 1\ar[r] &
  1\ar[r] &
    1&
        &
         &
          &
              \\
P_r\ar[u]\ar[r]^{\partial\rest{P_r}\ \ \ \ }&
    H^2({\coa{\sta C}},\pmmu_r)\ar[r]\ar[u]&
     H^2(\sta C, \pmmu_r)\ar[u]\ar[r]&
      P/rP\ar[u]\ar[r]&
        1 &
          &
           &
            \\
1\ar[r]\ar[u] &
    \Pic {\coa{\sta C}} \ar[u]\ar[r]^{\pi^*}&
     \Pic\sta C \ar[r]^{\sta j^*} \ar[u]&
      P\ar[r]\ar[u]&
       1\ar[r]\ar[u] &
         1 &
           &
              \\
1\ar[r]\ar[u] &
    \Pic {\coa{\sta C}} \ar[u]\ar[r]^{\pi^*}&
     \Pic\sta C \ar[r]^{\sta j^*}\ar[u]^r &
      P\ar[r]\ar[u]^r&
       1\ar[r]\ar[u] &
         1\ar[r]\ar[u] &
           1 &
              \\
1\ar[r]\ar[u] &
 (\Pic \coa{\stC})_r     \ar[r]\ar[u] &
  (\Pic \stC)_r     \ar[r]\ar[u] &
   P_r\ar[u]\ar[r]^{\partial\rest{P_r}\ \ \ \ }&
    H^2({\coa{\sta C}},\pmmu_r)\ar[r]\ar[u]&
     H^2(\sta C, \pmmu_r)\ar[u]\ar[r]&
      {\textstyle P/rP}\ar[u]\ar[r]&
        1\\
&
 1\ar[r]\ar[u] &
  1\ar[r]\ar[u] &
   1\ar[r]\ar[u] &
    \Pic {\coa{\sta C}} \ar[u]\ar[r]^{\pi^*}&
     \Pic\sta C \ar[r]^{\sta j^*}\ar[u]&
      P\ar[r]\ar[u]&
       1\ar[u]\\
&
 &
  1\ar[r]\ar[u] &
   1\ar[r]\ar[u] &
    \Pic {\coa{\sta C}} \ar[u]\ar[r]^{\pi^*}&
     \Pic\sta C \ar[r]^{\sta j^*}\ar[u]^r&
      P\ar[r]\ar[u]^{r}&
       \ar[u]\\
&
 &
   &
    1\ar[r]\ar[u] &
     (\Pic \coa{\stC})_r     \ar[r]\ar[u] &
      (\Pic \stC)_r     \ar[r]\ar[u] &
       P_r\ar[u]\ar[r]^{\partial\rest{P_r}}&
        \dots\ar[u]\\
&
 &
  &
    &
       1\ar[r]\ar[u] &
         1\ar[r]\ar[u] &
           1\ar[u] &
              \\
}}
\end{equation*}
\begin{caption}{the diagram is commutative and all horizontal and vertical sequences are exact.}\end{caption}
\end{figure}
}

In the figure above we used the following notation:
\begin{enumerate}
\item$P_{\sta e}$ denotes $\Pic(\sta B\Aut(\sta e))$, and
$P$ denotes ${\textstyle\prod_{\sta e} P_{\sta e}}$;
\item for any group $H$, $H_r$ denotes the $r$-torsion subgroup;
\item $\partial$ denotes the differential of
$\Ccal_\bullet(\Lambda, \ZZ/r\ZZ)$ for $\Lambda$ the dual graph of $\stC$ with an orientation;
\item $\sta j$ denotes the embedding of the singular locus in $\sta C$;
\item $\pi$ denotes the morphism from $\sta C$ to $\coa{\stC}$.
\end{enumerate}

We prove Theorem \ref{thm:rootsnum} in two steps:
we first focus on roots of $\Ocal$; then, we
consider the roots of $\sta F$.

\vspace{0.2cm}

\noindent{\it Step 1. The special case: $r$th roots of $\Ocal$.}
We need to show that
there are exactly $r^{2g}$ roots of $\Ocal$ if and only if $l(\sta e)$
is a multiple of $r$ for any nonseparating node $\sta e$.
By Theorem \ref{thm:roots} (see Figure 3.2), the number of elements of
$\Pic(\stC)_r$ equals the product of $\size{(\Pic \coa{\stC})_r}=r^{2g-1+\size{V}-\size{E}}$ by the
size of the kernel of
\begin{equation}\label{eq:P}
\partial \colon P_r \to H^2(\coa{\stC},\pmmu_r).
\end{equation}
Therefore, it is enough
to show the following claim
\begin{equation}\label{eq:claim}
\size{\ker\partial}=r^{1-\size{V}+\size{E}}
\quad  \Leftrightarrow \quad  l({\sta e})\in r\ZZ
\quad \forall \sta e \text{ nonseparating}.
\end{equation}
Recall that the
first Betti number of the
dual graph $\Lambda$ of $\stC$ is given by $b_1(\Lambda)=1-\size{V}+\size{E}$.
Consider the subgraph $\Lambda(\sta e)$ of
$\Lambda$ of ${\coa{\sta C}}$, whose vertices and edges are
$V$ and  $E\setminus \{\sta e\}$,
denote by $\widehat{\partial}_\sta e$ the restriction of
the chain differential of $\Lambda(\sta e)$ to
${P}_r({\sta e})=\prod_{E\setminus \{\sta e\}}(P_{\sta e})_r$.
For any nonseparating $\sta e\in E$, we have
\begin{equation}\label{eq:sizeker}\size{\ker \partial} =
\size{\ker \widehat{\partial}_{\sta e}}\cdot
\hcf \{r,l(\sta e)\}.
\end{equation}
The claim \eqref{eq:claim} follows.
Indeed, assume $\size{\ker\partial}=r^{1-\size{V}+\size{E}}$. Now, we have
$\size{\ker \widehat \partial_{\sta e} }\le
r^{b_1(\Lambda(\sta e))}$ and $b_1(\Lambda(\sta e))= b_1(\Lambda)-1$
if $\sta e$ is nonseparating.
Hence, $r=\hcf\{r,l(\sta e)\}$.
Conversely, assume $l({\sta e})\in r\ZZ$ for all nonseparating $\sta e$.
The claim holds when there are
$n-1$ nonseparating edges, then
\eqref{eq:sizeker} implies
$\size{\ker \partial}=r^{b_1(\Lambda(e))+1}=r^{1-\size{V}+\size{E}}$.

\vspace{0.2cm}

\noindent{\it Step 2. The general case: $r$th roots of  $\sta F$.}
We assume the numerical condition in the statement, and
we show that
it implies that $\sta F$
has
one $r$th root (by Step 1, this also implies that
the number of roots is $r^{2g}$).

By the diagram in Figure 3.2, this amounts to a combinatorial
criterion on the differential of $\Ccal^\bullet(\Lambda,\ZZ/r\ZZ)$
(recall that an orientation for $\Lambda$ is chosen).
The point is that
for any line bundle $A$ on the curve ${\coa{\sta C}}$
the pullback
$\pi^*A$ has an $r$th root in $\Pic(\stC)$
if and only if $\deg (A)$ is in the image of $\partial$
in Figure 3.2.

In order to state the criterion
we attach a partition of $E\setminus \{\sta e\}$ and $V$
to each
separating node $\sta e$ joining $v_+$ and $v_-$:
$$E\setminus \{\sta e\}=E^+\sqcup E^- \quad\quad
V=V^+\sqcup V^-, $$
where the set $E^+$
(the set $V^+$) contains the edges (the vertices) that can be
connected to $v_+$ without passing through $\sta e$.
Then, regard $P_r$ as the product $\prod_{\sta e \in E} (P_{\sta e})_r$
and write $P^\pm({\sta e})=\prod_{\sta e\in E^\pm}(P_{\sta e})_r$,
write $\partial_{\sta e} ^\pm$ for the corresponding chain differentials,
and denote by $\ep_{\sta e}^\pm$ the composition
of $V^\pm\into V$
with the augmentation homomorphism $\ep \colon V \to \ZZ/r\ZZ$
(over $\ker \varepsilon$, we obviously have $ \ep_{\sta e}^+=-\ep_{\sta e}^-$).
\begin{lem}\label{lem:comb}
Assume that $r$ divides $l(\sta e)$ if $\sta e$ is nonseparating.
Consider an element
$\lvec{t}\in \ker \varepsilon \in (\ZZ/r\ZZ)^V$;
the following conditions are equivalent.
\begin{enumerate}
\item \label{eq:inP}
$\lvec{t}$ is in the image of $P_r\subseteq (\ZZ/r\ZZ)^E$ via $\partial$.
\item \label{eq:aug} For any separating edge $\sta e$,
the value of $\ep_{\sta e}^+ (\lvec{t}\,)\in \ZZ/r\ZZ$ belongs to  $P_{\sta e}$
(or, equivalently, the value of $\ep_{\sta e}^- (\lvec{t})=-\ep_{\sta e}^+ (\lvec{t}\,)$
belongs to $P_{\sta e}$).
\end{enumerate}
\end{lem}
\begin{proof}
Note that the claim is trivial if
all edges are nonseparating. Indeed, \eqref{eq:aug}
is true, whereas
\eqref{eq:inP} is true, because
if $\ep(\lvec{t})$ vanishes,
then $\lvec{t}$ is exact in $\Ccal^\bullet$ and lies in $\partial (P_r)$ by
$l(\sta e)\in r\ZZ$.

Then, choose a separating edge $\sta e$ joining $v^+$ to $v^-$
Note that, the homomorphism
$\partial$
can be written as
\begin{align*}
P^+({\sta e})&\times (P_{\sta e})_r\times P^+({\sta e})
 &&\to && (\ZZ/r\ZZ)^{V^+}\quad\quad \times \quad \quad
(\ZZ/r\ZZ)^{V^-}\\
 &\ (a,x,b)&& \mapsto&& \big(\partial^+_{\sta e}(a)+ i_+(x),
\ \  \partial^-_{\sta e}(b)- i_-(x)\big),
\end{align*}
where $i^+$ and $i^-$ are the injections of
$(P_{\sta e})_r$ in $(\ZZ/r\ZZ)^{V^+}$
and $(\ZZ/r\ZZ)^{V^-}$ induced by $v^+\in V^+$ and
$v^-\in V^-$.

We assume $\lvec{t}\in \partial (P_r)$, and
we prove \eqref{eq:aug}. Indeed,
there exist $a\in P^+({\sta e})$  and $x\in (P_{\sta e})_r$
satisfying $\ep^+_{\sta e}(\lvec{t})=\ep^+_{\sta e}
(\partial^+_{\sta e}(a)+ i_+(x))$. Therefore, we have
$\ep^+_{\sta e}(\lvec{t})=\ep^+_{\sta e}
(i_+(x))=x\in P_{\sta e}.$

Conversely, choose a separating edge $\sta e$, and
set $x_{\sta e}=\ep_{\sta e}^+(\lvec{t})\in (P_\sta e)_r$.
We construct $\lvec{h}\in P_r$ such that $\partial (\lvec{h})=\lvec{t}$, by means of
the presentation of $\partial$ given in the diagram above.
Using the partition
$V=V^+\sqcup V^-$,
write $\lvec{t}$ as $(\lvec{t}^+_{\sta e},\lvec{t}^-_{\sta e})$.
Now, giving $\lvec{h}$ is equivalent to finding two elements mapping via
 $\partial_{\sta e}^+$ and  $\partial_{\sta e}^-$
to $\lvec{q}^+_{\sta e}=\lvec{t}^+_{\sta e}-i^+_{\sta e}(x_{\sta e})$ and
$\lvec{q}^-_{\sta e}=\lvec{t}^-_{\sta e}+i^+_{\sta e}(x_{\sta e})$.
By induction on the number of separating edges, the lemma holds on the subgraphs
with edges $E^+$ (or $E^-$)
and vertices $V^+$ (or $V^-$).
Therefore, in order to lift
$\lvec{q}=\lvec{q}^+_{\sta e}$ and $\lvec{q}^-_{\sta e}$
to $P^+(\sta e)$ and $P^-(\sta e)$,
we only need to show that $\lvec{q}$ satisfies \eqref{eq:aug}.
This is immediate, since, by construction, for any
separating edge $\sta f$
we have either $\ep_{\sta f}^+(\lvec{q})=
\ep_{\sta f}^+(\lvec{t})$
or $\ep_{\sta f}^-(\lvec{q})=\ep_{\sta f}^-(\lvec{t})$.
\end{proof}

Assuming that the numerical condition on the nodes is satisfied,
we construct an $r$th root of $\sta F$.
Using the orientation
of the dual graph $\Lambda$ chosen above, we
have a canonical way to associate to a node
$\sta e$ a point $\sta p_+$ over $\sta e$
lying in the normalization
of $\stC$ at $\sta e$. We set
$m(\sta e)=\mult_{\sta p_+}(\sta F)$.
In the same way, for any separating node $\sta e$, the orientation
allows us to associate to $\sta e$ a component $\sta C_+(\sta e)$ of the
partial normalization of $\stC$ at $\sta e$.
Then, we set $d(\sta e)=\deg_{\sta C_+(\sta e)}(\sta F)$.
By hypothesis there  exists
a function $k(\sta e)\in \ZZ$ satisfying
\begin{align}\label{eq:ke}
&d(\sta e)l(\sta e)=k(\sta e )r
&&\text{for separating nodes and}\\
&m(\sta e)=k(\sta e)r&&\text {for nonseparating nodes.}\label{eq:kebis}
\end{align}
The orientation of the dual graph $\Lambda$ also
induces a canonical generator
$\sta T_{\sta p_+}$
for each group $\Pic(\sta B\Aut(\sta e))$.
We consider
$$\textstyle{\prod_{\sta e}\sta T_{\sta p_+}^{\otimes k(\sta e)}
\in \prod_{\sta e} P_{\sta e}.}$$
Let $\sta M$ be a line bundle on $\stC$
such that
$\sta j^*\sta M=\textstyle{\prod_{\sta e}
\sta T_{\sta p_+}^{\otimes k(\sta e)}}$.
By the equations \eqref{eq:ke} and \eqref{eq:kebis} and Proposition \ref{pro:deg}, we have
$\sta j^*(\sta M^{\otimes r})=\sta j^*\sta F$.
Therefore, there exists a line bundle $A\in \Pic {\coa{\sta C}}$
satisfying $\pi^*A=(\sta M^{\otimes r})^\vee \otimes \sta F$.

In fact, $\pi^*A$ has an $r$th root, because
it lies in the kernel of $\Pic(\stC)\to H^2(\stC,\pmmu_r)$.
To see this, using the diagram of Figure 3.2,
it is enough to show that the
homomorphism $\Pic{\coa{\sta C}}\to H^2({\coa{\sta C}},\pmmu_r)$ sends $A$
into $\partial (P_r)$.
We apply Lemma \ref{lem:comb} to the multidegree of $A$ mod $r$;
clearly, condition (2) of the lemma holds if
the following
numerical condition is satisfied for each separating node:
the total degree of $\pi^*A$ on
the connected component $\wt{\stC}={\stC}_+(\sta e)$
is a multiple of $r$. Indeed, $r$ divides
\begin{multline*}
\deg(\pi^*A\rest{\wt{\stC}})=\deg((\sta M^{\otimes r})^\vee
\otimes \sta F)\rest{\wt{\stC}}
=-r\deg(\sta M\rest{\wt{\stC}}) +d(\sta e)
=r(-\deg(\sta M\rest{\wt{\stC}}) +d(\sta e)/r),
\end{multline*}
because  $\deg (\sta M\rest{\wt{\stC}})- k(\sta e)/l(\sta e)$
is an integer
by Proposition \ref{pro:deg}
and we have $k(\sta e)/l(\sta e)=d(\sta e)/r$ by \eqref{eq:ke}.
This proves the claim,
because by tensoring a root of $\pi^*A$
by $\sta M$ we get a root of $\sta F$.

Conversely, if $\sta F$ has $r^{2g}$ $r$th roots,
we show $l(\sta e), m(\sta e)\in r\ZZ$
for any nonseparating node $\sta e$ and
$l(\sta e)d(\sta e)\in r\ZZ$
for any separating node $\sta e$.
First,
for any nonseparating node $\sta e$,
Step 1 implies $l(\sta e)\in r\ZZ$,
and, since  $\sta F$ has an $r$th root,
$\sta j^*\sta F$ also has an $r$th root, and we have $m(\sta e)\in r\ZZ$.
Second, we show $l(\sta e)d(\sta e)\in r\ZZ$
for a separating node $\sta e$.
Since $\sta j^*\sta F$ has an $r$th root
and $\sta j^*$ is surjective,
we can choose a line bundle $\sta M$
on $\stC$ such that
$\sta j^*\sta M^{\otimes r}=\sta j^*\sta F$.
There exists $A\in \Pic{\coa{\sta C}}$ such that
$\pi^*A=(\sta M^{\otimes r})^\vee \otimes \sta F$.
Since $\sta F$ has an $r$th root,
$\pi^*A$ has an $r$th root; therefore,
the homomorphism $\Pic{\coa{\sta C}}\to H^2({\coa{\sta C}},\pmmu_r)$ sends $A$
into the image of
$\partial\colon P_r\to H^2({\coa{\sta C}},\pmmu_r).$
Note that,
by Lemma \ref{lem:comb},
the degree of $\pi^*A$ on $\stC_+(\sta e)$
has order $l(\sta e)$
modulo $r$; hence, we have
\begin{align*}
l(\sta e)\deg(\pi^*A\rest{\stC_+(\sta e)})\in r\ZZ,
\end{align*}
which implies $l(\sta e)d(\sta e)\in r\ZZ$, because
\begin{multline*}l(\sta e)\deg(\pi^*A\rest{\stC_+(\sta e)})
=l(\sta e)
\deg(\sta j^*(\sta M^{\otimes r})^\vee \otimes \sta F)\rest{\stC_1}
=-rl(\sta e)\deg(\sta M\rest{\stC_1}) +l(\sta e)d(\sta e).
\end{multline*}
\qed

\section{The notion of stability for twisted curves}
\label{sect:olsson}
\subsection{Twisted curves and the notion of $\lvec{l}$-stability}
By Olsson's Theorem \ref{thm:Olsson},
the category of twisted curves of genus $g\ge 2$ forms a Deligne--Mumford stack.
As the following example shows,
the stack $\wt \MMM_g$ is not separated.
\begin{exa}
A twisted curve $\stC$ over a discrete valuation ring $R$
with smooth generic fibre $\stC_K$ is isomorphic to
its coarse space over $K$ and may differ from it on the special fibre;
in this case  the coarse space ${\coa{\sta C}}$ and $\stC$  are two nonisomorphic
twisted curves extending $\stC_K$ on $S$.
Therefore, the valuative criterion of separateness
fails.
\end{exa}

Inside $\wt\MMM_g$, we identify all proper substacks containing $\MMM_g$.
We need to recall
the following standard notion of type of a node.
\begin{notn}[type of a node]\label{notn:type}
Given a twisted curve of genus $g\ge 2$ over an algebraically closed field
$k$ and a node $e\in \stC$
we set the following convention:
$$
e\in \stC \text{ is } \begin{cases}
\text{of type $0$ \quad}&\text{if the normalization
of $\stC$ at $\sta e$}\\
& \text{is connected (i.e. $\sta e $ is nonseparating);}\\
\\
\text{of type $i$ with \quad\quad\quad}& \text{if
normalizing  $\stC$ at $\sta e$ we get $\stC_1\sqcup \stC_2$
}\\
\text{$1\le i\le\lfloor g/2\rfloor$}& \text{with
$\{g(\stC_1),g(\stC_2)\}=\{i,g-i\}.$\quad}\end{cases}$$
\end{notn}
\begin{defn}\label{defn:lstab}
Let $\lvec{l}=(l_0,l_1,\dots,l_{\lfloor g/2\rfloor})$ be a
multiindex of positive and invertible integers $l_i$.
A twisted curve $\stC\to X$
is \emph{$\lvec{l}$-stable}, if the coarse space is stable and if
the stabilizer at a node of type $i$ has order $l_i$.
\end{defn}
\begin{theor}\label{thm:olsson}
Let us denote by $\MMM_g(\lvec l)$ the category
of $\lvec{l}$-stable curves. It is contained in
$\wt \MMM_g$ and it contains $\MMM_g$:
$$\MMM_g \into \MMM_g(\lvec{l})\into \wt \MMM_g.$$

\begin{enumerate}
\item[I.] The stack $\MMM_g(\lvec{l})$ is  tame, proper
(separated), smooth, irreducible and of Deligne--Mumford type.
The morphism $\MMM_g(\lvec{l})\to \MMMbar_g$ is finite, flat,
and is an isomorphism on the open dense substack $\MMM_g$.
\item[II.] Any proper substack $\sta X$ of $\widetilde
\MMM_g$ fitting in $\MMM_g \into \sta X\into \wt \MMM_g$ is
isomorphic to $\MMM_g(\lvec{l})$ for a suitable multiindex
$\lvec{l}$.
\end{enumerate}
\end{theor}
\begin{proof}
There is a natural
surjective morphism
$$\wt \MMM_g\to \MMMbar_g,$$
which is the functor sending
a twisted curve $\stC \to X$ to the coarse space $\coa{\stC}\to X$.
Note that $\MMM_g$ is dense in $\wt{\MMM}_g$: by
Theorem \ref{thm:Olsson}
any twisted curve
$\stC\to \Spec k$ can be realized
as the special fibre of a twisted curve
$\stC'$ over a discrete valuation ring in such a way that
the generic fibre is smooth.

Point (I) follows
from \cite[Thm.~1.9]{Ol} and, in particular, from the description of
versal deformation spaces \eqref{eq:versal}.
The morphism $\MMM_g(\lvec{l})\to \MMMbar_g$
is locally represented by the flat, finite, tame
morphism of Deligne--Mumford type
\begin{equation}\label{eq:localolsson}
[(\Spec \wt I)/\pmmu_{h_1}\times \dots \times \pmmu_{h_m}]\to \Spec I
\end{equation}
where $I$ is the versal deformation space of
a point of $\MMMbar_g$, $m$ is the number of nodes
of the curve represented by such point,
$\wt I=I[z_1, \dots, z_m]/(z_1^{h_1}
 - t_1, \dots , z_m^{h_m} - t_m)$ and
$\pmmu_{h_1}\times\dots \times  \pmmu_{h_m}$ acts
as $(\xi_{h_1},\dots,\xi_{h_m})z_i=\xi_{h_j}z_j$
(note that the index $h_j$ depends on the type of the
 $j$th node and
$h_j=l_i$ if
the $j$th node is of type $i$).
This means that $\MMM_g(\lvec{l})$ is smooth, which is also
shown in \cite[\S3]{ACV} and \cite[\S3]{AJ}.
The fact that $\MMM_g(\lvec{l})$ is irreducible is a consequence of the
fact that $\MMM_g$ is dense on $\MMM_g(\lvec{l})$.

We show (II). Let $\sta X$ be a proper stack, which contains
$\MMM_g$ and is contained in $\wt \MMM_g$.
By restriction of $\wt \MMM_g\to \MMMbar_g$ to $\sta X$, we obtain
a morphism of proper stacks
$\sta f\colon\sta X\to \MMMbar_g$.
Note that $\MMM_g$ is dense in $\sta X$ and in $\MMMbar_g$.
The valuative criterion of properness for $\sta X$ and
$\MMMbar_g$ implies that
for any geometric point $\sta y\colon\Spec k\to\MMMbar_g$
there exists a point $\sta x\colon \Spec k\to \sta X$
lifting $\sta y$. In fact $\sta x$ is unique.
To see this, consider a versal deformation $\Spec I \to \MMMbar_g$
at $\sta y$. The base change of $\sta X\to \MMMbar_g$ via
$\Spec I\to \MMMbar_g$ is a restriction of \eqref{eq:localolsson}.
The point $\sta y$ is in the locus $z_1=z_2=\dots=z_m=0$ and admits only
one lifting.

Now, we define positive indexes $\lvec{l}=(l_0,l_1,\dots,l_{\lfloor
g/2\rfloor})$ such that $\sta X$ is isomorphic to $\MMM_g(\lvec l)$. Denote by
$\sta u\colon \sta U({\MMMbar_g})\to \MMMbar_g$ the universal
stable curve. Consider a geometric point $\sta p\colon \Spec k \to \sta
U({\MMMbar_g})$, which is  a node of type $i$. Let $\sta q\colon
\Spec k\to \sta X$ be the unique morphism
which  lifts $\sta u\circ \sta p$.
The object determined by $\sta q$ is a twisted curve $\stC$,
whose coarse space $\coa{\stC}$ is represented by the point $\sta u\circ \sta p$.
We define $l_i$ as the order of the stabilizer of $\stC$ on the
node $\sta p$. The index $l_i$
is locally constant on the (connected) substack of $\sta
U({\MMMbar_g})$ of nodes of type $i$. Therefore, $l_i$ only
depends on $i$. This implies that
the objets of $\sta X$ are $\lvec{l}$-stable curves.
Finally, the properness of $\sta X$ implies
that any point $\Spec k\to \MMM_g(\lvec{l})$ lifts to $\sta X$.
Since $\MMM_g(\lvec l)$ is regular,
this suffices to show that $\sta X$ is isomorphic to $\MMM_g(\lvec l)$.
\end{proof}

In fact, the stack $\MMM_g(\lvec{l})$
admits an  alternative description.

The first part of Olsson's proof of Theorem \ref{thm:Olsson} in \cite{Ol}
consists of constructing a functor from
the category of twisted curves
to the category of
simple logarithmic extensions
of the logarithmic structure $\mathcal M_\Delta$ on $\MMMbar_g$
canonically
associated to the boundary locus $\Delta=\MMMbar_g\setminus \MMM_g$.
Indeed, the boundary locus is a
normal crossings divisor $\Delta=\sum _i\Delta_i$
where $\Delta_i$ is
the full subcategory of stable curves $C\to X$ where the
geometric fibres on every $x$ in $X$ are curves containing at least
a node of type $i$ (recall that a logarithmic structure is
canonically associated to any normal crossings divisor
\cite[1.5]{Ka}).
The second part of Olsson's proof shows that
simple extensions of logarithmic structures form an algebraic stack.

We point out that,  by  \cite[Lem.~5,3]{Ol}, Olsson's functor
can be regarded as an
equivalence between  the \emph{subcategory} $\MMM_g(\lvec{l})$ of
$\lvec{l}=(l_0,l_1,\dots, l_{\lfloor g/2\rfloor})$-stable curves
and the \emph{subcategory} of
simple extensions $\mathcal M_{\Delta}\to \mathcal N$
of locally free logarithmic structures
inducing, at each geometric point $x\colon \Spec k\to X$,
the commutative diagram
\[\xymatrix@C=2cm{
\NN^{\oplus k}\ar[d]\ar[r]^{\oplus_{i=1}^k \times a_i}& \NN^{\oplus k}\ar[d]\\
\ol{\mathcal M}_{\Delta,x}\ar[r]&\ol{\mathcal N}_x
}\]
where $\{1, \dots ,k\}$ is the set of
irreducible components of the
local picture of $\Delta$ at $x$, and
$a_i=l_j$ if $i\in \{1,\dots, k\}$
corresponds to the component $\Delta_j$ in $\Delta$
(following Olsson
we adopt the notation $\ol{\mathcal M}_{\Delta}=\mathcal M_{\Delta}/\Ocal^{\times }$,
$\ol{\mathcal N}=\mathcal N/\Ocal^{\times }$
and we refer to \cite[Lem~4.2]{MO} for the proof of the
canonical decomposition
of $\ol{\mathcal M}_{\Delta,x}\cong \NN^{\oplus k}$).

Now, we notice that
this subcategory of logarithmic extensions  is
precisely the category used in
Matsuki and Olsson's Generalization \ref{gen:log}.
In this way, we get the following statement.
\begin{theor}
For any multiindex $(l_0,l_1,\dots,l_{\lfloor g/2\rfloor})$,
we have the following isomorphism
\[\xymatrix@C=0cm{
\MMM_g(\lvec{l})\ar[dr]&\cong &\MMMbar_g(\textstyle{\sum_i \Delta_i/l_i})\ar[dl],\\
&\ \ \MMMbar_g&
}\]
where
$\MMM_g(\lvec{l})$ is the stack of $\lvec{l}$-stable curves
and
$\MMMbar_g(\sum_i \Delta_i/l_i)$
is the stack of simple extensions
of the logarithmic structure associated to
$\sum_i \Delta_i$
with indexes
$l_i$ in the sense of \cite[(4.2.2)]{MO}. \qed
\end{theor}

\subsection{The moduli stack of $r$th roots on $\lvec{l}$-stable curves}\label{sect:global}
Let ${\sta {LB}}_g$ be
the category
formed by pairs $(C\to X, M)$,
where  $C\to X$
is a smooth curve and
$M$ is a line bundle on $C$.
The stack  ${\sta {LB}}_g$ is the category fibred over
$\MMM_g$ whose fibre over a smooth curve
$f\colon C\to X$ is the stack $\sta{LB}_f$.
Let $\sta F$ be a line bundle on the universal curve
of $\MMM_g$, whose relative degree is a multiple of $r$.
Let $\MMM_g^{\sta F, r}$ be the fibred category whose fibre over
$f\colon C\to X$ is ${\sta F^{1/r}}$.
By Proposition \ref{pro:stackLBF/r},
$\MMM_g^{\sta F, r}$ is a Deligne--Mumford stack, \'etale over $\MMM_g$.

It is well known that $\sta F$ is
a power $\omega^{\otimes k}$
of the relative dualizing sheaf
modulo pullbacks from $\MMM_g$ (Enriques and
Franchetta's conjecture, \cite{Ha} \cite{Me} \cite{AC}).
In view of a compactification of $\MMM_g^{\sta F, r}$, we focus on
the case
$$\sta F=\omega^{\otimes k}$$
and we assume
$$(2g-2)k\in r\ZZ.$$

We extend $\sta {LB}_g$ to $\wt\MMM_g$:
we consider the category
$\wt{\sta {LB}}_g$ formed by pairs $(\stC\to X, \sta M)$
where $\stC\to X$
is a twisted curve and $\sta M$ is a line bundle on $\stC$, Remark \ref{rem:globalLB}.
We get the stack
\begin{equation}\label{eq:LBg}
\wt{\sta {LB}}_g\to \wt\MMM_g.
\end{equation}
By abuse of notation, we denote by $\sta F$
the power of the relative dualizing sheaf $\omega^{\otimes k}$
on the universal curve over $\wt\MMM_g$.
Let $\wt \MMM{}_g^{\sta F, r}$ be the fibred category over
$\wt \MMM_g$ whose fibre over $\sta f\colon \stC\to X$ is the stack
of ${\sta F^{1/r}}$ of $r$th roots of $\sta F$
on $\sta C$.
In this way, the objects of $\wt \MMM_g^{\sta F, r}$
are triples $(\stC\to X,\stL,\sta j)$,
where $\stC$ is a twisted curve, $\stL$ is a line bundle on
$\stC$ and $\sta j$ is an isomorphism $\sta j\colon \stL^{\otimes r}
\to \sta F_{\stC}$. Morphisms from $(\sta C\to X,\stL, \sta j)$ to
 $(\sta C'\to X',\stL', \sta j')$
are pairs
$(\sta m,\sta a)$ where $\sta m\colon \sta C\to \stC'$
is a morphism of twisted curves and
$\sta a\colon \stL\to \sta m^*\stL'$
is an isomorphism of line bundles, with $\sta a^{\otimes r}$
commuting with $\sta j$ and $\sta j'$.

By Proposition \ref{pro:stackLBF/r}, $\wt \MMM_g^{\sta F, r}$ is a
Deligne--Mumford stack \'etale over $\wt \MMM_g$
\begin{equation}\label{eq:LBF/r}
\wt \MMM_g^{\sta F, r}\to \wt \MMM_g,
\end{equation}
which can be regarded as the projection on the second factor
of the fibred product
$(\wt{\sta {LB}}_g)\,{} _{{\sta k}_r}\!\!\times _{\sta F} \wt \MMM_g$
fitting in
\[\xymatrix@R=0cm{
\wt \MMM_g^{\sta F, r}\ar[rr]\ar[dd]&& \wt{\sta {LB}}_g\ar[dd]^{\sta k_r}\\
&\square&\\
\wt\MMM_g\ar[rr]_{\sta F}&& \wt{\sta{LB}}_g
}\]
where $\sta k_r$ is induced by
the $r$th power in $\GG_m$ and
$\sta F=\omega^{\otimes k}$
is regarded as
a section of \eqref{eq:LBg}.

We consider the \'etale morphism
$\wt\MMM_g^{\sta F, r}\to \wt \MMM_g$.
Its restriction to the stack $\MMM_g^{\sta F, r}$
of $r$th roots of $\sta F$ on smooth curves $\sta C\to X$
forms a finite stack
on $\MMM_g$ equipped with a
torsor structure under
the group stack $\MMM_g^{\Ocal, r}$.
On the other hand $\wt \MMM_g$ is not separated and
so is $\wt\MMM_g^{\sta F, r}$.
By Theorem \ref{thm:olsson}, we can consider the base changes
to all compactifications of $\MMM_g$ in $\wt \MMM_g$
via $\MMM_g(\lvec{l})\into \wt \MMM_g$:
we get
$\MMM_g^{\sta F, r}(\lvec{l})\to \wt \MMM_g(\lvec{l})$
\[\xymatrix@R=0.1cm{
\MMM_g^{\sta F, r}\ar[dd]\ar[rr] &&
\MMM_g^{\sta F, r}(\lvec{l})\ar[dd]\ar[rr]&& \wt \MMM_g^{\sta F, r}\ar[dd]\\
&\square& &\square&\\
\MMM_g\ar[rr] && \MMM_g(\lvec{l})\ar[rr] &&\wt\MMM_g.
}\]
We now characterize the compactifications
$\MMM_g^{\sta F, r}(\lvec{l})$ for which the
morphism to $\wt \MMM_g(\lvec{l})$
is proper.
 \begin{theor}\label{thm:COND}
For any $\sta F=\omega^{\otimes k}$,
the category  $\MMM^{\sta F, r}_g(\lvec{l})$ is a smooth
Deligne--Mumford algebraic stack, \'etale on $\MMM_g(\lvec{l})$.
\begin{enumerate}
\item[I.] For $\sta F=\Ocal$,
the stack $\MMM_g^{\Ocal , r}(\lvec{l})$ is a finite group stack
if and only if $r$ divides $l_0$.
\item[II.] For $\sta F=\omega$ and $2g-2\in r\ZZ$,
the stack $\MMM_g^{\Ocal , r}(\lvec{l})$ is a finite
group stack and $\MMM_g^{\omega , r}(\lvec{l})$ is a finite  torsor
under $\MMM_g^{\Ocal , r}(\lvec{l})$ if and only if $r$ divides
$$(2i-1)l_i \text { for all $i$}.$$
In this way, we obtain several compactifications
of the stack $\MMM_g^{\omega, r}$ of smooth $r$-spin curves:
for each $\lvec{l}$ satisfying $l_i(2i-1)\in r\ZZ$,
$$\MMM_g^{\omega, r}(\lvec{l})\to \MMM_g(\lvec{l})$$
is the finite torsor of \emph{$r$-spin $\lvec{l}$-stable curves}.
\item[III.] More generally, for $\sta F=\omega^{\otimes k}$  and $(2g-2)k\in r\ZZ$,
the stack $\MMM_g^{\Ocal , r}(\lvec{l})$ is a finite group stack
and $\MMM_g^{\sta F , r}(\lvec{l})$ is a finite  torsor
under $\MMM_g^{\Ocal , r}(\lvec{l})$ if and only if
$r$ divides
$$l_0 \quad \text{and} \quad (2i-1)kl_i, \text { for $i>0$}.$$
\end{enumerate}
\end{theor}
\begin{proof}
It is enough to show (III). We check the numerical condition of
Theorem \ref{thm:rootsnum} for $\sta F=\omega^{\otimes k}$.
Let $\sta e$ be a node of a twisted curve.
If $\sta e$ is  nonseparating, the condition amounts to
require that $r$ divides $\size{ \Aut(\sta e)}$.
If $\sta e$ is separating and of type $i$,
then the numerical condition is that $r$ divides $\size {\Aut(\sta e)}k(2i-1)$.
In this way, Theorem \ref{thm:rootsnum}
implies the claim.
\end{proof}
\begin{pro}\label{pro:torsors}
The category $\MMM^{\Ocal, r}_g(\lvec{l})$
is equivalent to the category of $\pmmu_r$-torsors on
$\lvec{l}$-stable curves. In this way, as soon as $r$ divides $l_0$,
the stack  $\MMM^{\Ocal, r}_g(\lvec{l})$ is a compactification
of the stack of  $\pmmu_r$-torsors on smooth curves.
\end{pro}
\begin{proof}
There is a natural  functor from $\pmmu_r$-torsors to
$r$-torsion line bundles. Given a torsor $\sta T$ on $\sta C$
with an action of  $\pmmu_r$, consider
$\sta P=\sta T\times _\stC \Af^1$ over $\sta C$.
Note that $\pmmu_r$ acts on both factors and the
diagonal action on the fibre product is free. The quotient
$\sta P/\pmmu_r$ yields an $r$-torsion line bundle on $\stC$.
The functor is essentially surjective and fully faithful.
\end{proof}
\begin{rem}\label{rem:globalrigi}
For each object of $\MMM_g^{\sta F, r}(\lvec{l})$ there is an
injection of $\pmmu_r$ in the  automorphism group
(the $r$th roots of unity
act by multiplication along the fibres of the line bundle).
The rigidification of  $\MMM_g^{\sta F, r}(\lvec{l})$
along the group scheme $\pmmu_r$
yields a representable cover $\sta H({\sta F})^{\pmmu_r}$
of $\MMM_g(\lvec{l})$, see
Proposition \ref{pro:rigidified}.
In this way,
the morphism $\MMM^{\sta F, r}_g(\lvec{l})\to \MMM_g(\lvec{l})$
factors as
$$\MMM^{\sta F, r}_g(\lvec{l})\to \sta H({\sta F})^{\pmmu_r}
\to  \MMM_g(\lvec{l}),$$
where $\sta H({\sta F})^{\pmmu_r}\to\MMM_g(\lvec{l}) $
is a representable cover of degree $r^{2g}$ and
$\MMM^{\sta F, r}_g(\lvec{l})\to \sta H({\sta F})^{\pmmu_r} $
is an \'etale $\pmmu_r$-gerbe.
By Proposition \ref{pro:rigidified} and Proposition
\ref{pro:finiteness}, as soon as
$\lvec{l}$ satisfies
$l_0\in r\ZZ$ and $(2i-1)kl_i\in r\ZZ$ for $i>0$
the base change of
$\sta H({\Ocal})^{\pmmu_r}\to\MMM_g(\lvec{l}) $
and of $\sta H({\sta F})^{\pmmu_r}\to\MMM_g(\lvec{l}) $
with respect to a morphism from a scheme $X$ to $\MMM_g(\lvec{l})$
yields a finite group $X$-scheme $G_X$ and
a finite torsor $T_X$ under $G_X$.
\end{rem}
\begin{rem}[$\la$-stable curves]
Note that Theorem \ref{thm:COND} is automatically satisfied if
$l_i=\la$ for all $i$ and $\la$ is a multiple of $r$.
In this case, we write ``$\la$-stable'' instead
$\lvec{l}$-stable.
\end{rem}
\begin{rem}[$n$-pointed curves]
Olsson's results \cite[Thm.~1.9]{Ol}
are formulated in terms of $n$-pointed curves;
therefore, the category $\wt \MMM_{g,n}$
of twisted curves, with $n$ ordered distinct markings
$\sigma_1,\dots, \sigma_n$ in the smooth
locus, is a Deligne--Mumford stack.
We point out that our method produces proper
stacks of $r$th roots of any line bundle $\sta F=\omega^{\otimes k}(\sum_i -h_i[\sigma_i])$
on the universal twisted curve over $\wt{\MMM}_{g,n}$
for all integers $k, h_1, \dots, h_n$ satisfying $(2g-2)k-\sum_i h_i\in r\ZZ$.
(The case $k=1$ is relevant to Witten's conjecture \cite{Wi}).
This happens by Corollary \ref{cor:rdividesl} and the fact that
$\sta F=\omega^{\otimes k}$ is a pullback from the
universal (scheme-theoretic) stable $n$-pointed curve over $\wt{\MMM}_{g,n}$,
Proposition \ref{pro:omega}.
In general, Corollary \ref{cor:rdividesl} yields the following statement.
\end{rem}
\begin{cor}\label{cor:F/rmarked}
Let $\sta F$ be a line bundle on the universal twisted  $n$-pointed curve
isomorphic to a pullback from the universal stable  $n$-pointed curve.
Let $\wt \MMM_{g,n}^{\sta F, r}$ be
the fibred category on $\wt\MMM_{g,n}$
of $r$th roots of $\sta F$ on twisted curves.
It is a stack, \'etale on $\wt \MMM_{g,n}$.
It is nonempty, as long as we assume that the relative degree of $\sta F$ is a multiple of $r$.

For any $\la\in \ZZ$, the stack of $\la$-stable
$n$-pointed curves $\MMM_{g,n}(\la)$
is smooth, irreducible, and proper.
If $r$ divides $\la$,
the stack $\MMM^{\Ocal, r}_{g,n}(\la)$
of $r$-torsion line bundles on $\la$-stable
$n$-pointed curves
is a finite group stack
and the stack
$\MMM_{g,n}^{\sta F, r}(\la)$
of $r$th roots of $\sta F$ on $\la$-stable curves is a finite  torsor
under the group stack $\MMM_{g,n}^{\Ocal , r}(\la)$.

In this way, we obtain several compactifications
of the stack $\MMM_{g,n}^{\omega, r}$ of smooth $r$-spin $n$-pointed curves:
for each $\lambda\in r\ZZ$,
$$\MMM_{g,n}^{\omega, r}(\lambda)\to \MMM_{g,n}(\lambda)$$
is the finite torsor of \emph{$r$-spin $n$-pointed $\lambda$-stable curves}.
\qed
\end{cor}

\begin{rem}
Let $r_1$ and $r_2$ be positive and coprime integers.
In \cite[\S4, Rem.~4.11]{JKV2}, it
is noted that the functor
$L\mapsto (L^{\otimes r_2}, L^{\otimes r_1})$
is an isomorphism
$\MMM^{\sta F,r_1r_2}_{g}\cong \MMM^{\sta F, r_1}_{g}\times_{\MMM_{g}} \MMM^{\sta F, r_2}_{g}$, that does not
extend to the
compactifications given in the existing literature. With our formalism the
isomorphism extends immediately (we omit markings for simplicity).\end{rem}
\begin{pro}[roots of two coprime orders]\label{pro:r1r2}
Let $r_1$ and $r_2$ be positive and coprime integers. Set
$\la=r_1r_2$. The functor
$\stL\mapsto (\stL^{\otimes r_2}, \stL^{\otimes r_1})$
is an isomorphism of stacks
$$ \MMM^{\sta F,r_1r_2}_g (\la)\cong \MMM^{\sta F,r_1}_g(\la)\times_{\MMM_g(\la)}\MMM^{\sta F,r_2}_g(\la).$$
\end{pro}
\begin{proof}
The inverse functor is $(\stL_1,\stL_2)\mapsto \stL_1^{\otimes h_2}
\otimes \stL_2^{\otimes h_1}$ for $h_1$ and $h_2$ satisfying
$h_1r_1+h_2r_2=1$. Indeed, we have
$(\stL_1^{\otimes h_2} \otimes \stL_2^{\otimes h_1})^{\otimes r_1r_2}=
\stL_1^{\otimes r_1r_2h_2} \otimes \stL_2^{\otimes
r_2r_1h_1}\cong\sta F^{\otimes
h_1r_1+h_2r_2}= \sta F.$
\end{proof}

\begin{exa}
[$g=1$ and $n=1$]\label{exa:roots}
Consider the twisted curve $\stC=[E/\pmmu_2]$, where
$E$ is equal to $\PP^1/(0\sim \infty)$ and
$\pmmu_2$ acts by change of sign.
We now exhibit the four distinct
square roots of $\omega_\stC$ up to isomorphism.
If we denote by $\sta x$ a smooth point of $\sta C$,
we can regard this example as
a check that the fibre of
${\MMM}^{\omega/2}_{1,1}(2)\to {\MMM}_{1,1}(2)$ over the
geometric point representing $(\stC,\sta x)$
actually contains $4$ distinct geometric points.

Note that $\omega_\stC$ is trivial, so we are actually looking for
square roots of $\mathcal O_{\stC}$.
Consider the normalization $\stC^\nu\to \stC$, which is
isomorphic to $[\PP^1/\pmmu_2]$. On $\stC^\nu$ there are $2$ roots of
$\mathcal O$: the line bundle
$\mathcal O_{\PP^1}$ with trivial $\pmmu_2$-action
on the fibres; and
the line bundle $\mathcal O_{\PP^1}$ with
$\pmmu_2$-action on the fibres given by $t\mapsto -t$.
Each of these line bundles
descend to $E$ to form a root of $\mathcal O_E$
in exactly two nonisomorphic ways.
Therefore, on $\stC$, we get four square roots $\stL^{++}$, $\stL^{+-}$, $\stL^{-+}$, and $\stL^{--}$
of $\omega_\stC$ up to isomorphism.

Note that the line bundles $L^{\sigma,\tau}$ above,
are $2$-torsion line bundles and can be regarded as
the $2$-torsion subgroups of $\Pic^{\pmmu_2}(E)$, the
group of  $\pmmu_2$-linearized line bundles.
Their geometric realizations are stacks fibred over $\stC$, which we describe explicitly.
Let
$W^+$ be the line bundle on $E$,
obtained from $\PP^1\times \Af^1$ by glueing the lines over $0$ and $\infty$ via
$(0,t)\sim(\infty,t)$
and let $W^-$ be the line bundle on $E$
obtained from $\PP^1\times \Af^1$ by glueing the lines over $0$ and $\infty$ via
$(0,t)\sim(\infty,-t)$.
The geometric bundles over $\stC$
associated to $\sta L^{\sigma,\tau}$
for ${\sigma,\tau}\in \{+,-\}$
is the quotient stack $[W^{\tau}/\pmmu_2]$ with $\pmmu_2$
acting as $t\mapsto {\sigma} t$ on the fibres.

Since $\omega_\stC$ is isomorphic to $\Ocal$, following Proposition \ref{pro:torsors},
we also provide a concrete
description of the line bundles above in terms of
$\pmmu_2$-torsors on $\stC$. We exhibit four distinct
representable $2$-folded \'etale covers of $\stC$.
Note that some points of the covering stack $\sta D$ over
$\stC$ might have
nontrivial stabilizer $\pmmu_2$:
in the drawings we adopt the convention of
marking them with a black circle. The labelling
$\sta D^{++}$, $\sta D ^{+-}$, $\sta D^{-+}$, and $\sta D^{--}$
matches the above notations for the line bundles.

\smallskip
\noindent\emph{1. The cover $\sta D^{++}$.} It is
$\sta D^{++}=\stC\times \pmmu_2 $ on which $\pmmu_2$ acts
as $\id\times \xi_2$.
The morphism $\sta p$ is the projection
to the first factor.

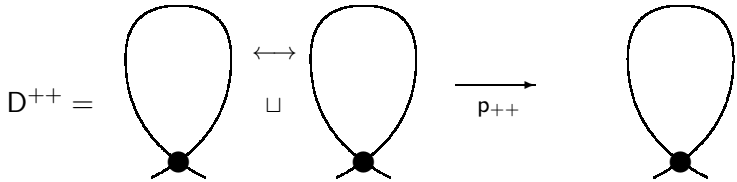
\begin{picture}(200,80)(-50,0)
  \qbezier(15,5)(45,20)(45,50)
  \qbezier(35,5)(5,20)(5,50)
  \qbezier(5,50)(5,70)(25,70)
  \qbezier(25,70)(45,70)(45,50)
  \put(25,11){\circle*{8}}
  \put(95,11){\circle*{8}}
  \put(215,11){\circle*{8}}
  \put(58,30){\footnotesize
$\sqcup$}
  \put(53,50){\footnotesize
$\longleftrightarrow$}
  \put(-40,30){
$\sta D^{++}=$}
  \qbezier(85,5)(115,20)(115,50)
  \qbezier(105,5)(75,20)(75,50)
  \qbezier(75,50)(75,70)(95,70)
  \qbezier(95,70)(115,70)(115,50)
  \put(138,30){\footnotesize
$\sta p_{++}$}
  \put(130,40){\vector(1,0){30}}
  \qbezier(205,5)(235,20)(235,50)
  \qbezier(225,5)(195,20)(195,50)
  \qbezier(195,50)(195,70)(215,70)
  \qbezier(215,70)(235,70)(235,50)
 \end{picture}

\smallskip
\noindent\emph{2. The cover $\sta D^{+-}$.} Take
$[\PP^1/\pm]\times \pmmu_2 $ modulo the relation
$(0,\sigma)\sim (\infty,-\sigma)$, for all $\sigma\in \pmmu_2$.
The $\pmmu_2$-action is generated by
$\id\times \xi_2$ and
the morphism $\sta p$
is the projection to the first factor.

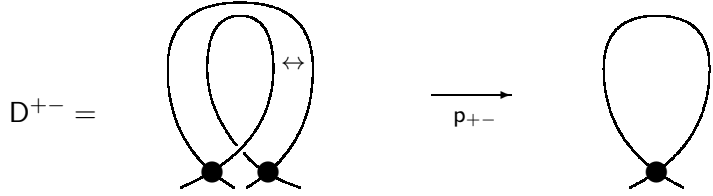
\begin{picture}(200,80)(-50,0)
  \qbezier(35,5)(70,20)(70,50)
  \qbezier(55,5)(30,20)(30,50)
  \qbezier(30,50)(30,75)(57.5,75)
  \qbezier(45,50)(45,70)(57.5,70)
  \put(47,11){\circle*{8}}
  \put(68,11){\circle*{8}}
  \put(215,11){\circle*{8}}
  \put(73,50){\footnotesize
$\leftrightarrow$}
  \put(-30,30){
$\sta D^{+-}=$}
  \qbezier(60,5)(85,20)(85,50)
  \qbezier(80,5)(65,10)(59,18)
  \qbezier(56,21)(45,35)(45,50)
  \qbezier(57.5,70)(70,70)(70,50)
  \qbezier(57.5,75)(85,75)(85,50)
  \put(138,30){\footnotesize
$\sta p_{+-}$}
  \put(130,40){\vector(1,0){30}}
  \qbezier(205,5)(235,20)(235,50)
  \qbezier(225,5)(195,20)(195,50)
  \qbezier(195,50)(195,70)(215,70)
  \qbezier(215,70)(235,70)(235,50)
 \end{picture}

\smallskip
\noindent\emph{3. The cover $\sta D^{-+}$.} Take
the \'etale atlas $E=\PP^1/(0\sim \infty)$ of $\stC$.
The  $\pmmu_2$-action is the change of sign on $\PP^1$.
The morphism $\sta p$ is $E\to [E/\pmmu_2]$.
Note that, the local picture of $\sta p$
at the node is given by
$(z,w)\mapsto (z,w)$ on
$\{zw=0\}\to [\{z'w'=0\}/\pmmu_2]$.

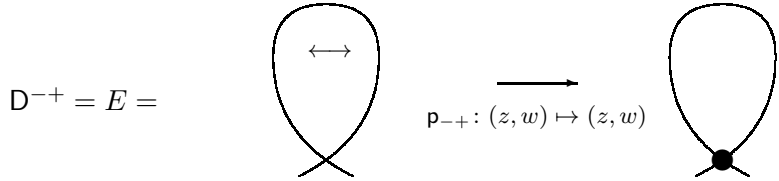
\begin{picture}(200,80)(-50,0)
  \qbezier(55,5)(85,20)(85,50)
  \qbezier(75,5)(45,20)(45,50)
  \qbezier(45,50)(45,70)(65,70)
  \qbezier(65,70)(85,70)(85,50)
  \put(130,40){\vector(1,0){30}}
  \put(215,11){\circle*{8}}
  \put(58,50){\footnotesize
$\longleftrightarrow$}
  \put(-55,30){
$\sta D^{-+}=E=$}
  \put(103,25){\footnotesize $\sta p_{-+}\colon (z,w)\mapsto(z,w)$}
  \qbezier(205,5)(235,20)(235,50)
  \qbezier(225,5)(195,20)(195,50)
  \qbezier(195,50)(195,70)(215,70)
  \qbezier(215,70)(235,70)(235,50)
 \end{picture}

\smallskip
\noindent\emph{4. The cover $\sta D^{--}$.} Take
$E=\PP^1/(0\sim \infty)$ with $\pmmu_2$
acting by change of sign as above.
On the smooth locus
$E\sm\longrightarrow \stC\sm$
the morphism $\sta p$
is $x\mapsto x^2$.
On the other hand, we define $\sta p$ so that its local picture
at the node  $\{zw=0\}\to [\{z'w'=0\}/\pmmu_2]$
is given by
$(z,w)\mapsto (z,-w)$.

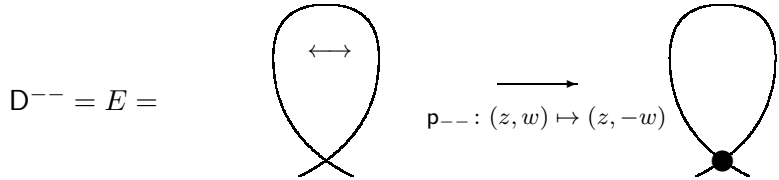
\begin{picture}(200,80)(-50,0)
  \qbezier(55,5)(85,20)(85,50)
  \qbezier(75,5)(45,20)(45,50)
  \qbezier(45,50)(45,70)(65,70)
  \qbezier(65,70)(85,70)(85,50)
  \put(130,40){\vector(1,0){30}}
  \put(215,11){\circle*{8}}
  \put(58,50){\footnotesize
$\longleftrightarrow$}
  \put(-55,30){
$\sta D^{--}=E=$}
  \put(103,25){\footnotesize $\sta p_{--}\colon (z,w)\mapsto(z,-w)$}
  \qbezier(205,5)(235,20)(235,50)
  \qbezier(225,5)(195,20)(195,50)
  \qbezier(195,50)(195,70)(215,70)
  \qbezier(215,70)(235,70)(235,50)
 \end{picture}

This completes the check that the fibre
of $\MMM_{1,1}^{\omega/2}(2)\to \MMM_{1,1}(2)$
over the curve $(\stC,\sta x)$ in
$\MMM_{1,1}(2)$ is the $0$-dimensional stack
given by $4$ disjoint copies of $\BBB(\ZZ/2\ZZ)$
(recall that each object has a nontrivial automorphism
acting by multiplication by $-1$ along the fibre of the
line bundle).

One can ask a natural question at this point: what is
the fibre of the corresponding morphism between coarse spaces
over the closed point corresponding
to $(\stC,\sta x)$ in the moduli space
of stable curves?
In order to answer this question one should
note that the only automorphism of $(\stC, \sta x)$
that acts nontrivially on the objects of $\MMM_{1,1}^{\omega/2}(2)$
is the automorphism $\sta g$ of order $2$ generating $\Aut(\stC,\coa{\sta C})$.
Then, by Proposition \ref{pro:aut_lb}, we note that the action
of $\sta g$ fixes  $\sta L^{++}$ and $\sta L^{+-}$ and
identifies the objects $\sta L^{-+}$ and $\sta L^{--}$.
Indeed the formula in Proposition \ref{pro:aut_lb} can be written
as
\begin{align*}\sta g^*\stL^{++}=\stL^{++}\otimes \stL^{++}=\stL^{++},\\
\sta g^*\stL^{+-}=\stL^{+-}\otimes \stL^{++}=\stL^{+-},\\
\sta g^*\stL^{-+}=\stL^{-+}\otimes \stL^{+-}=\stL^{--}.
\end{align*}
Therefore, the fibre is the disjoint union of
two reduced point and a third point of length two over
the coarse moduli space of stable curves.
\end{exa}
\begin{rem}
By Proposition \ref{pro:aut_lb}, this analysis can be generalized
to any integer $r>2$. For simplicity, we consider $r$ prime.
Let $\stC$ be the $r$-stable curve on the $1$-pointed
nodal curve $\PP^1/(0\sim \infty)$.
The group $\Aut(\stC,\coa{\stC})=\pmmu_r$ acts
freely on the $(r^2-r)$ spin structures that are not pullbacks from $\coa{\stC}$.
This means that the fibre of the morphism to $\ol \MMM_{1,1}$
contains $1$ point representing the trivial $r$-spin structures and
$2r-2$ points representing nontrivial $r$-spin structures, half of which are pullbacks
from the coarse space.
After the identifications induced by the hyperelliptic involution we get
\begin{align}\label{eq:counter}
&\text{there are exactly $r-1$ nontrivial
$r$-spin $r$-stable curves over $\PP^1/(0\sim \infty)$.}
\end{align}
This allows to picture the coarse space of
nontrivial $r$-spin curves of genus $1$, with $1$ marking.
This leads us to point out the following counterexample to
Conjecture 4.2.1 of \cite{Ja_Pic}, which predicts that the Picard group of
the stack of smooth $r$-spin $1$-pointed curves of genus $1$ is finite.
\end{rem}
\begin{exa} \label{exa:counter} The space
$N_r$ of nontrivial $r$-stable $r$-spin
curves is a curve covering $\coa{\MMM_{1,1}}$:
indeed it is an $(r^2-1)/2$-fold cover of the projective line.
For instance, fix $r=11$; then $N_r\to \coa{\MMM_{1,1}}$ has degree 60.
Over the two curves with extra automorphisms there are
respectively $(r^2-1)/4=30$ and $(r^2-1)/6=20$ spin structures.
By \eqref{eq:counter}, there are exactly $r-1=10$ singular spin curves.
Then, by the Riemann--Hurwitz formula, the Euler--Poincar\'e characteristic
is
$\chi(N_r)=0.$
In this way, $N_r$ is a genus-$1$ curve. The moduli stack of
nontrivial $r$-spin structures on smooth $1$-pointed genus-$1$ curves
is a stack over the genus-$1$ curve $N_r$ minus a finite number of points.
Its Picard group cannot be finite.
A similar computation using \eqref{eq:counter}
shows that for any prime integer
$r>3$ we have
\begin{equation}
g(N_r)=(r-5)(r-7)/24.
\end{equation}
\end{exa}

\subsection{The relation  with Abramovich and Jarvis's compactification
}\label{sect:preexist}
The compactifications \cite{Ja_geom} and \cite{AJ} adopt two
different methods but are isomorphic, \cite[Prop.~4.3.1]{AJ}). We
restate the construction. We use
systematically the equivalence between line bundles on a stack
$\sta X$ and morphisms $\sta X\to \BBB\GG_m$.

The compactification $\sta B_{g,n}(\sta
B\GG_m,\omega_{log}^{1/r})$ introduced in \cite{AJ} is
the following category. An {\em object} is the datum of a
$1$-commutative diagram
\[ \xymatrix@C=1cm{
                    & \stC\ar[dl]_{\stL}
                    \ar[dr]^{\omega_{\stC/X}(\sum_{i} \sta S_i)}&  \\
 \BBB\GG_m\ar[rr]_{\sta k_r}&          &  \BBB \GG_m \\
} \] where $\sta k_r$ is induced by the homomorphism $t\mapsto
t^r$ and the following conditions are satisfied.
\begin{enumerate}
\item The stack $\stC$ is of Deligne--Mumford type,
flat of relative dimension $1$ with nodal singularities over $X$.
\item The stacks $\sta S_1, \dots, \sta S_n$ are
closed disjoint substacks of $\stC\sm$ and \'etale gerbes over
$X$.
\item The corresponding coarse spaces $\coa{\stC}, \coa{\sta S_1},\dots, \coa{\sta S_n}$
form a proper, $n$-pointed, nodal curve over $X$, and
$\pi \colon
\stC\to \coa{\stC}$ is an isomorphism away from the nodes and the stacks
$\sta S_i$.
\item At a node $\sta p$ in $\stC$,
for a suitable integer $l$, the local picture is given by
$[V/\pmmu_{l}]$, where, for some $t\in T$, $V$ is
$\Spec(T[z,w]/(zw - t))$ and $\pmmu_l$ acts as
$(z,w)\mapsto (\xi_lz,\xi_l^{-1}w)$.
\item The morphism $\stL$ is representable.
\end{enumerate}
A {\em morphism} is a $1$-commutative diagram
\[ \xymatrix@C=1cm{
& \stC\ar[d]\ar[ddl]_{\stL}
                    \ar[ddr]&  \\
                    & \stC'\ar[dl]^{\stL'}
                    \ar[dr]&  \\
 \BBB\GG_m\ar[rr]_{\sta k_r}&          &  \BBB \GG_m .\\
} \]
As usual, morphism are considered up to
$2$-isomorphisms (see Lemma \ref{lem:AV}).

In fact, we show that $\sta B_{g,n}(\sta
B\GG_m,\omega_{log}^{1/r})$ is a compactification of
$$\bigsqcup _{0\le h_i <r } \MMM_{g,n}^{\omega(\pmb h), r},$$
where $\omega(\pmb h):=\omega(-\sum_i h_iS_i)$
is the relative dualizing sheaf on the universal curve
over $\MMM_{g,n}$ twisted by the divisors
$S_n$ determining
the $i$th marking.
\begin{defn}[faithful line bundles]
A line bundle  $\sta M$
on a twisted curve $\sta C\to X$ is faithful if
it satisfies the following condition: for each node
$\sta e$ on $\stC$ the action of
$\Aut(\sta e)$ on $\sta M_{\sta e}$ is faithful.
\end{defn}
\begin{rem}
In the case of a line bundle $\sta M$ whose $r$th tensor power
is the pullback of a line bundle $F$
on the coarse space, for any
$g\in \Aut(\sta e )$ the element $g^r$
acts as the identity on $\sta M_{\sta e}$.
Therefore, the fact that $\sta M$ is faithful implies that
$\size{\Aut(\sta e)}$ divides $r$.
\end{rem}
\begin{pro}\label{pro:compAJ}
There is a finite and surjective morphism
\begin{equation}
\label{eq:mortoAJ}
\bigsqcup _{0\le h_i <r }
\MMM_{g,n}^{\omega(\pmb h), r}(r)
\longrightarrow \sta B_{g,n}(\sta
B\GG_m,\omega_{log}^{1/r}),\end{equation}
where
$\MMM_{g,n}^{\omega(\pmb h), r}(r)$
is the stack of $r$th roots of $\omega(\pmb h)
=\omega(-\sum_i h_iS_i)$
over $n$-pointed $r$-stable curves.

The morphism (\ref{eq:mortoAJ}) has degree one, but, in general, it is not an isomorphism.
Indeed, consider a point
$\sta x$ in $\MMM_{g,n}^{\omega(\pmb h), r}(r)$ and
its image $\sta y$ in $\sta B_{g,n}(\sta
B\GG_m,\omega_{log}^{1/r})$.
We have
$$\size{\Aut(\sta x)}/
\size{\Aut(\sta y)}=r^m/\textstyle{\prod_{i=1}^m
d_i}$$
where $d_1,\dots , d_m$
are the orders of
the automorphism groups
of
the nodes $\sta e_1,\dots,\sta e_m$ of the twisted curve corresponding to $\sta y$.
\end{pro}
\begin{proof}
There is an equivalence of categories
$$\sta B_{g,n}(\sta
B\GG_m,\omega_{log}^{1/r})
\longrightarrow  \bigsqcup _{0\le h_i <r }
\sta P(h_1,\dots, h_n),$$
where $\sta P(h_1,\dots, h_n)$ is the category on
$\wt \MMM_{g,n}$ of faithful $r$th roots of
the line bundle $\omega(\pmb h)=\omega(-\sum_i h_iS_i)$  on twisted curves.
We prove the equivalence.

First, we need to introduce a decomposition of $\sta B_{g,n}(\sta B\GG_m,\omega_{log}^{1/r})$
into connected components ${\sta P}({\lvec l},{\lvec m})$.
The local picture of a line bundle
$\sta L$ on $\stC$ at a point $\sta p\colon \Spec k\to \sta
S_i\into \stC$ is given, for suitable indexes $l_i$ and $m_i$, by
a $\pmmu_{l_i}$-equivariant line bundle  $\sta W$ on
$[V/\pmmu_{l_i}]$ where $V$ is $\Spec T[z]$, the action is
$z\mapsto \xi_{l_i}z$, and $\sta W$ is linearized by the
character $\xi_{l_i}\mapsto \xi_{l_i}^{m_i}$. The coefficients
$l_i$ and $m_i$ induce locally constant morphisms $\sta l,\sta
m\colon\sta B_{g,n}(\sta B\GG_m,\omega_{log}^{1/r})\to
\ZZ^{n}$ with $0\le \sta m< \sta l$. Note also that the
representability condition on $\sta L$ implies that $\sta
l$ and $\sta m$ are coprime on each coordinate of $\ZZ^n$. The
values taken by $\sta l$ and $\sta m$ determine a  decomposition
of $\sta B_{g,n}(\sta B\GG_m,\omega_{log}^{1/r})$ as the
disjoint union of substacks ${\sta P}({\lvec l},{\lvec m})$ where
the vectors $\lvec l=(l_i), \lvec m=(m_i)\in \ZZ^n$ satisfy $0\le
\lvec l< \lvec m$ and $\hcf \{m_i,l_i\}=1$ for any $i$.

Second, there is an equivalence of category between ${\sta P}({\lvec l},{\lvec m})$
and the stack  ${\sta P}(h_1,\dots,h_n)$
for $h_i=rm_i/l_i-1$.
By \cite[Thm.~4.1]{Ca} and \cite[Thm.~1.8]{Ol}, there is an equivalence
between the category of twisted curves equipped with
$n$ distinct sections
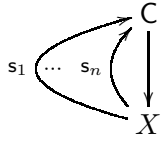
\[ \xymatrix@R=1cm{
{\stC}\ar[d]\\
X\ar@/^1.5cm/[u]^{\sta s_1}_\cdots\ar@/^0.5cm/[u]^{\sta s_n}\\
} \]
in the smooth locus
and the category of Deligne--Mumford stacks $\stC\to X$
satisfying conditions (1-4) in the definition of
$\sta B_{g,n}(\sta B\GG_m,\omega_{log}^{1/r})$,
where we require that the automorphism groups
have order $l_i$ over the points of $S_i$.
The functor
sends $(\stC\to X,\sta s_1,\dots,\sta s_2)$ to
the stack $\sta C\times_ {\coa{\sta C}} \sta D$, where
we set $\sta D={\coa{\sta C}}[S_1/l_1]\times_{\coa{\sta C}} \dots \times _{\coa{\sta C}} {\coa{\sta C}}[S_n/l_n]$
for $S_i=\sta s_i(X)$. Note that
the latter stack is equipped with a natural
projection to $\stC$. For $h_i=rm_i/l_i-1$,
we have a functor
${\sta P}(h_1,\dots,h_n)\to {\sta P}({\lvec l},{\lvec m}) $
induced by pushforward via the projection. Indeed,
the pushforward of line bundles on $\sta D$
is a line bundle on ${\coa{\sta C}}$, see \cite[\S3]{AJ} or \cite{Ca}.
In this way, the functor lands in
${\sta P}(h_1,\dots,h_n)$ (the representability assumption at the nodes
is equivalent to requiring that
the $r$th root is faithful at the nodes).
The inverse functor is induced by pullback
and tensorization with the tautological line bundles
$\sta M_i^{\otimes m_i}$, which are defined on $\sta D$
and satisfy $\sta M_i^{\otimes l_i }=\Ocal (S_i)$.

The morphism \eqref{eq:mortoAJ}
is the composite of  ${\sta P}(h_1,\dots,h_n)\to {\sta P}({\lvec l},{\lvec m})$ and
the disjoint union of finite and surjective morphisms
of stacks
\begin{equation}\label{eq:AJrestricted}
\MMM_{g,n}^{\omega(\pmb h), r}(r)\to \sta P(h_1,\dots h_n),
\end{equation}
whose restriction to the open and dense substack
$\MMM_{g,n}^{\omega(\pmb h), r}$
is the identity and whose corresponding morphism
between coarse spaces is an isomorphism.
We define the morphisms \eqref{eq:AJrestricted} here below.

The functor
$\MMM_{g,n}^{\omega(\pmb h), r}(r)\to \sta P(h_1,\dots h_n)$,
sends the morphism of stacks $\stC\xrightarrow{L}
\BBB\GG_m$ to the corresponding representable morphism
$\stC'\xrightarrow{L'} \BBB\GG_m$ making
$\stC\to \stC'\to \BBB\GG_m$ the
the ``relative moduli
space''
in the sense of \cite[5.2.4,(c)]{ACV}.
By means of the weak
valuative criterion it is easy to see that this functor defines a
surjection (each geometric point is lifted as in the  proof of Theorem \ref{thm:olsson}).
The functor sending the object determined by
$\stC$ and $\stL$ to the object determined by
$\stC'$ and $\stL'$ is not an isomorphism in general. Indeed,
the ratio of the orders of the automorphism groups at the two objects
is $\size {\Aut(\stC,\stC')}$, where
$\Aut(\stC,\stC')$ denotes the group of automorphisms of $\stC$
that fix $\stC'$ (this happens because $\sta L$ is the pullback of
$\sta L'$ via $\sta C\to \sta C'$).
By Theorem \ref{thm:aut_tw}, the order of $\Aut(\stC,\stC')$ is equal to
$r^m/\textstyle{\prod_{i=1}^m
d_i}$, if we denote by
$d_1,\dots , d_m$
the orders of
the automorphism groups
of
the nodes $\sta e_1,\dots,\sta e_m$ of $\stC'$.
\end{proof}
\begin{rem}[the Witten top Chern class]
The Witten top Chern class is a rational Chow cohomology class
which plays a crucial role in the definition of the relevant
numerical invariants in Witten's conjecture \cite{Wi}. Although
the new compactifications are not isomorphic to the preexisting
one, the surjective morphism above yields an isomorphism between
the coarse spaces.
This implies that the rational Chow rings are
isomorphic, see for example \cite{Kr}.

 There are two equivalent
formulations \cite{PV} and \cite{Ch}
of the construction of the
Witten top Chern class
and they both use the universal stable
$r$-spin structure of Jarvis's compactification \cite{Ja_geom},
which is a sheaf of rank $1$ rather than an invertible sheaf.
Nevertheless they can be applied without modification to the new
compactifications and yield the same class
after the identification of the rational Chow cohomology rings.
\end{rem}
\begin{pro}\label{pro:wtcc}
The Witten top Chern class functor defined in \cite{PV} and \cite{Ch}
yields a class $c_{\rm W}$ in the rational cohomology
of $\MMM_{g,n}^{\omega(\pmb h), r}(r)$ as
well as a class $\ol c_{\rm W}$
in the rational cohomology of Abramovich and Jarvis's compactification.
The outputs are compatible in the sense that
$c_{\rm W}$ is a pullback of $\ol c_{\rm W}$ via the
surjective morphism
of degree one exhibited in Proposition \ref{pro:compAJ}.
\end{pro}
\begin{proof}
To see this, note that both
constructions start from a datum in the derived category which is
obtained by pushing Jarvis's universal
$r$th root $\mathcal L$ and the
universal homomorphism $f$ along
the universal stable curve. The
morphism between our compactification and the compactification of
Abramovich and Jarvis induces a morphism between the universal
twisted curve on $\MMM^{\omega(\pmb h), r}_{g,n}(r)$ and the universal
stable curve on Jarvis's compactification. We only need to check that the pushforward to the
universal stable curve yields the
universal sheaf-theoretic stable $r$-spin structure
$(\mathcal L,f)$ of Jarvis's construction \cite{Ja_geom}.
The proof of this fact can be found in
\cite[\S3, \S4.3]{AJ}, which applies verbatim to our setting.
\end{proof}

\appendix
\section{Appendix. The stack $\sta {LB}_{\sta f}$}
We analyse the category
$\sta{LB}_{\sta f}$ of line bundles
on a flat and proper morphism $\sta f\colon \sta Y\to X$,
where $\sta Y$ is a tame stack of Deligne--Mumford type,
and $\coa{\sta f}\colon \coa{\sta Y}\to X$
is flat on $X$.
The fact that this category $\sta{LB}_{\sta f}$ forms a stack even
when $\sta f$ is not represented by a scheme is
a preliminary to our compactification.
We show that it can be proven by adapting
Mumford's treatment of cohomology and base change \cite[II.~\S5]{Mu} to
the stack-theoretic situation: $\sta f\colon \sta Y\to X$.
As mentioned in \S3.1, M.~Lieblich provides a
more general statement implying
this result by
showing that
the category of flat families of coherent sheaves on $\sta Y$
with proper support over $X$ is an algebraic stack,
\cite[Thm.~2.1.1, Lem.~2.3.1]{Li}.

\subsection{The fibred category
$\sta {LB}_{\sta f}$ is a stack}
Let $\sta Y$ be a
Deligne--Mumford stack, flat and proper
on a base scheme $X$
$$\sta f\colon \sta Y\to X.$$
We write
$\sta{LB}_{\sta f}$
for the category of line bundles on
base changes $\sta Y_S=\sta Y\times_X S$
for every $X$-scheme $S$.
More precisely, the \emph{objects} are pairs $(S,\sta M)$, where
$S$ is an $X$-scheme
and $\sta M$ is a line bundle on
$\sta Y_S=\sta Y\times_X S$. The
\emph{morphisms} $(S,\sta M)\to (S',\sta M')$
are pairs $(m,\sta a)$, where $m\in \mathsf{Hom}_X(S,S')$
and $\sta a$ is an isomorphism of line bundles
$\sta a\colon \sta M\xrightarrow{\sim} \sta M'\times _{S'} S$
on $\sta Y_S$.
\begin{rem}
The category $\sta{LB}_{\sta f}$
induces a functor
sending an  $X$-scheme $S$ to the \emph{groupoid } in
$\sta{LB}_{\sta f}$ formed by the objects $(S,\sta M)$ on $S$.
Even when $\sta f$ is
a representable morphism $f\colon Y\to X$,
such a functor differs from the functor
sending an  $X$-scheme $S$ to
the \emph{set} $\Pic (Y_{S})$, which is
the functor used
in Grothendieck's
treatment \cite{FGA}
of the relative Picard functor.
We illustrate the relation between the two functors
at the end of this appendix using
the notion of rigidification of a stack along
a group scheme, Theorem \ref{thm:rigid}.
\end{rem}
\begin{pro}\label{thm:lb}
On a base scheme $X$, let
$\sta f\colon \sta Y\to X$ be
a flat and proper morphism of Deligne--Mumford type, with geometrically connected fibres, tame, and
coarsely represented by
an $X$-scheme $\coa{\sta Y}$, projective and flat on $X$.
Then, the category $\sta{LB}_{\sta f}$ is a stack on $X$.
\end{pro}
\begin{proof}
Since the category
is fibred in groupoids by definition, we only need to show
\begin{enumerate}
\item
the representability of the isomorphism functors,
\item
the effectiveness of any \'etale descent
datum of objects
\end{enumerate}
(as far as $\wt {\sta{LB}}_g$ is concerned, we ignore the
issue of geometrization, namely
the existence of a smooth and surjective
morphism from a scheme to
the stack, see Remark \ref{rem:algebraic}).
We show point (1) by means of
the following statement, a
stack-theoretic generalization
of Mumford's theorem on cohomology and base change
\cite[II.~\S5]{Mu} for schemes.
\begin{lem}\label{lem:cohbc}
Let $\sta Y$ be a tame
Deligne--Mumford stack, flat and proper over
an affine scheme $S$. Assume that the morphism of schemes
$\coa{\sta Y}\to S$ is flat.
Let $\Ecal$ be a locally free and coherent sheaf
on $\sta Y$.
Then, there exists a perfect complex
$K^\bullet: 0\to K^0\to K^1\to \dots \to K^n\to 0$ on $S$ and an isomorphism of
functors
$$H^p(\sta Y_T,\Ecal_T)=H^p(K^\bullet_T)$$
on the category of $S$-schemes $T$ (here,
$\sta Y_T$, $\Ecal_T$, and
$K^\bullet_T$ denote the base change via $T\to S$).
\end{lem}
\begin{proof}
Recall that the direct image via
$\pi\colon \sta Y\to \coa{\sta Y}$
is an exact functor from the category of coherent sheaves
on $\sta Y$ to the category of coherent sheaves on $\coa{\sta Y}$,
\cite[Lem.2.3.4]{AV}.
Therefore, we have
the isomorphism  $H^p(\sta Y_T,\Ecal_T)\cong H^p(\coa{\sta Y}_T,\pi_*\Ecal_T)$,
and it is enough to
find a complex $K^\bullet$
satisfying
$$H^p(\coa{\sta Y}_T,\pi_*\Ecal_T)\cong H^p(K^\bullet_T).$$
Mumford' theorem \cite[II.~\S5]{Mu}
shows that such a $K^\bullet$ exists
if $\pi_*\Ecal$ is
a coherent sheaf, flat on $S=\Spec A$.
Indeed, we check that $\pi_*\Ecal$ is flat on $S$, which
means that, locally on $\coa{\sta Y}$, there exists an affine open set
on which $\pi_*\Ecal$ is given by a flat $A$-module.
This happens
because, as shown in \cite[Lem.2.2.3]{AV},
for a stack  $\sta Y$ of Deligne--Mumford type,
there is an \'etale covering $Y_\al\to \coa {\sta Y}$
such that, for all $\al$, the
pullback $\sta Y\times _{\coa {\sta Y}} Y_\al $
is a quotient stack of the form $[U_\al/G_\al]$,
where $U_\al$ is a scheme and $G_\al$ is a finite group
acting on $U_\al$
(note that since $\sta Y$ is tame, $\size{G_\al}$ is prime to the
residue characteristic).
So, on an affine open set $V\subset \coa{\sta Y}$, $\pi_*\Ecal$
can be regarded as the direct image of a locally free
$G$-equivariant coherent sheaf
on an  affine scheme $U=\Spec R$.
Such a sheaf can be regarded as a $G$-linearized $R$-module $M$.
Therefore, as an $\Ocal_V$-module, $\pi_*\Ecal$ corresponds
to $M^G$, the submodule of $M$ of $G$-invariant elements.
The tameness assumption
implies linear reductiveness:
$M$ splits as  $M^G\oplus M'$.
So, as an $A$-module, $M^G$ is flat,
because it is a direct summand of $M$
and, on the other hand,
$M$ is flat over $R$, which is flat over $A$.
\end{proof}

Now, point (1) follows from the following lemma.
\begin{lem}\label{lem:isom}
For any scheme $S$ and objects $\al=(S,\sta M_\al)$ and
$\be=(S,\sta M_\be)$
in $\sta{LB}_{\sta f}$
the functor $\mathsf{Isom}_S(\al,\be)$
from $S$-schemes to sets is represented by a separated
scheme
locally of finite type over $S$.
\end{lem}
\begin{proof}
Write
$\sta D_S=\sta M_{\al}\otimes \sta M_{\be}^\vee$.
We write $\pi$ for  $\sta Y\to \coa{\sta Y}$, and
we adopt the notation
$\sta p\colon \sta Y_S\to S$ and
$p\colon \coa{\sta Y}_S\to S$.
For any morphism $Z\to S$ write $\sta Y_Z$ and $\sta D_Z$
for the base change of $\sta Y_S$ and $\sta D_S$.
We need to represent the functor $\Isom_S(\al,\be)$
sending an $S$-scheme $Z$
to the set of isomorphisms between $\sta M_\al$ and $\sta M_\be$
on $Z$. This is equivalent to the set
of nonzero sections $s\in\Gamma(\sta Y_Z,\sta D_Z)$,
where $s$ is nowhere vanishing.
By Lemma \ref{lem:cohbc},
we can take $d\colon K^0\to K^1$,
a homomorphism of vector bundles on $S$;
then, the $S$-scheme $T_{\al,\be}=\{d=0\}$ of $K^0$
represents represent the functor $\sta{Hom}_S(\al,\be)$
sending an $S$-scheme $Z$
to the set of homomorphisms between $\sta M_\al$ and $\sta M_\be$
on $Z$. The composition of homomorphism induces a morphism
$$c\colon T_{\al,\be}\times_S T_{\be,\al} \to T_{\al,\al}.$$
Let $1\colon S\to T_{\al,\al}$ be the section representing the identity
homomorphism.
The functor \linebreak $\Isom_S(\al,\be)$ is represented by $\{c=1\}$.
Since
$K^0$ is a finite dimensional vector bundle,
$T_{\al,\be}$
is separated and of finite type over $S$.
It follows that the scheme representing the
functor $\Isom_S(\al,\be)$ is separated
and of finite type (it is closed in the scheme $T_{\al,\be}\times_S T_{\be,\al}$,
which is separated and of finite type over $S$).
\end{proof}

Finally, we show point (2): any \'etale descent datum of line bundles on
$\sta Y$ is effective. Indeed, given an $X$-scheme
$S$, an \'etale cover $(S_\al\to S)$,
and objects $(S_\al, \sta M_\al)$, where $\sta M_\al$ is a line bundle
on $\sta Y_\al=\sta Y\times_S S_\al$,
together with isomorphisms between pullbacks of
$(S_\al, \sta M_\al)$ and $(S_\be, \sta M_\be)$ to
$S_\al\times_S S_\be$ satisfying
the cocycle condition, we claim that these data descend
to an object $(S,\sta M)$, where $\sta M$
is a line bundle on $\sta Y_S$.
Indeed $\sta M$ is defined by
\'etale descent of line bundles and morphisms of line bundles for schemes:
for any scheme $T\to \sta Y_S$ the line bundle $\sta M_T$ on $T$
is induced by descent along
the \'etale cover $(\sta Y_{\al}\times_{\sta Y_S} T\to T)$
of the line bundles
$\sta M_\al\times_{\sta Y_S} T$ on $\sta Y_{\al}\times_{\sta Y_S} T$.
\end{proof}

\begin{rem}\label{rem:algebraic}
Although we did not show that $\sta{LB}_{\sta f}$ is
an algebraic stack in the sense of Artin's definition, we point out that
Lemma \ref{lem:isom} together with the fact that the
relative cotangent complex of $\sta F^{1/r}\to X$
is trivial (see \ref{pro:stackLBF/r}) can be used to
prove that
the stack $\sta F^{1/r}$ is algebraic.  The
claim follows easily from Artin's method, in which one
starts from a deformation theory,
constructs formal
deformation spaces, and shows that they are algebraizable.
Indeed the deformation functor of $\sta F^{1/r}$ coincides with that
of $X$, since the relative cotangent complex of $\sta F^{1/r}\to X$
is trivial. To finish the proof one only needs to show
that formal deformations are algebraizable, which is an immediate
consequence of
 Grothendieck's Existence Theorem extended to tame stacks by Abramovich and Vistoli, \cite{AV}.
\end{rem}

\subsection{Rigidification}
Note that each object $(S,\stL)$ of $\sta{LB}_{\sta f}$
over an $X$-scheme $S$ has automorphisms given by multiplication
by $s\in \Gamma (S,\GG_m)$ along the fibre of $\stL$.
More precisely, $\sta H=\sta {LB}_{\sta f}$
and $G=\GG_m$
fit in the following setting.
\begin{theor}[Abramovich, Corti, Vistoli, {\cite[Thm.~5.1.5]{ACV}}]
\label{thm:rigid}
Let $\sta H$ be a stack on a base scheme $X$,
let $G$ be a flat finitely presented group scheme on $X$,
and assume that for any object $\tau$ of $\sta H$
over an $X$-scheme $S$
there is an embedding
$i_\tau\colon G(S)\into \Aut_S(\tau)$
compatible with pullbacks
in the obvious sense (for any $\fie\colon \tau\to \tau'$ in $\sta H$ over
the morphism of $X$-schemes $f\colon S\to S'$
we have $i_\tau \circ f^*=\fie^* \circ i_{\tau'}$).
Then, there exists a stack $\sta H^{G}$ and
a morphism of stacks $\sta H\to \sta H^G$
over $X$
satisfying the following conditions.
\begin{enumerate}
\item For any object $\tau \in \sta H (S)$
with image $\xi\in \sta H^G(S)$, the
set $G(S)$ lies in the kernel of
$\Aut_S (\xi) \to \Aut_S (\tau)$.
\item The morphism $\sta H\to \sta H^G$
is universal for morphisms of stacks $\sta H\to \sta H'$
satisfying (1) above.
\item In the condition (1) above, if $S$ is the spectrum of an
algebraically closed field, we have
$\Aut_{S}(\xi) = \Aut _{S}(\tau)/G(S)$.
\end{enumerate}
If $\sta H$ is an algebraic stack, then $\sta H^G$ is also
an algebraic stack. If $\sta H$ is of Deligne--Mumford type,
then $\sta H^G$ is also of Deligne--Mumford type and the
coarse space
$\coa{\sta H}$ is isomorphic to $\coa{\sta H^G}$.
\end{theor}
We call $\sta H^G$ the
\emph{rigidification of $\sta H$ along $G$}.
In \cite[I.~Prop.~3.0.2, (2)]{Ro}, Romagny shows that if
$\sta H$ is a Deligne--Mumford stack, then
it
is ``locally isomorphic'' to
$\sta BG$ (on $\sta H^G$)
and is indeed an \'etale  $G$-gerbe.

The construction of $\sta H^G$ consists of two steps.
\begin{enumerate}
\item We define a prestack ${\sta H}_{\rm pre}^G$: the
category whose objects are the objects of $\sta H$ and whose morphisms
are obtained by means of a quotient operation on
the sheaves of morphisms of $\sta H$ (for any object $\tau$
the embeddings
$G(S)\into \Aut_S(\tau)$ induce a
categorically injective morphism of $S$-group schemes
of the pullback $G_S$ of $G$ to $S$
to the group scheme $\Aut_S(\tau)$ of automorphisms of $\tau$).
\item We pass to the stack $\sta H^G$ associated to the prestack
${\sta H}_{\rm pre}^G$ in the
sense of \cite[Lem.~3.2]{LM}.
\end{enumerate}

This construction provides a natural framework to a standard procedure that
occurs systematically in the construction of the Picard functor.
In general, for any morphism of schemes $f\colon Y\to X$,
the natural functor $S\mapsto \Pic(Y_S)$ from $X$-schemes
to sets is a presheaf and is not represented by a scheme.
The actual ``relative Picard functor'' is defined
by the passage to
the associated sheaf.
This point is illustrated in detail in \cite[Ch.~8]{BLR}, by
Bosch, L\"utkebohmert, and Raynaud.
In this way, the construction of the
relative Picard functor is just another way to
rigidify $\sta{LB}_f$ along $\GG_m$.

\small
{

\end{document}